\documentclass[a4paper]{amsart}

\usepackage{amsmath,amssymb} 
\usepackage{comment}

\newtheorem{theorem}{Theorem}[section]
\newtheorem{lemma}[theorem]{Lemma}
\newtheorem{mycor}[theorem]{Corollary}

\theoremstyle{definition}
\newtheorem{definition}[theorem]{Definition}
\newtheorem{example}[theorem]{Example}
\newtheorem{xca}[theorem]{Exercise}

\theoremstyle{remark}
\newtheorem{remark}[theorem]{Remark}

\numberwithin{equation}{section}

\newcommand{\abs}[1]{\lvert#1\rvert}

\newcommand{\blankbox}[2]{%
  \parbox{\columnwidth}{\centering
    \setlength{\fboxsep}{0pt}%
    \fbox{\raisebox{0pt}[#2]{\hspace{#1}}}%
  }%
}

\newbox\squ  
\setbox\squ=\hbox{\vrule width.6pt
   \vbox{\hrule height.6pt width.4em\kern1ex
         \hrule height.6pt}%
   \vrule width.6pt}
\def\sqbox{\copy\squ}
\def\proof{\noindent{\sl Proof.\quad}}
\def\endproof{\quad\sqbox\par\vspace{2mm}}
\def\mod{\operatorname{-mod^{fd}}}
\def\id{\operatorname{id}}

\def\ZZ{{\mathbb Z}}
\def\fn{{\operatorname{fin}}}
\def\P{{\mathcal P}}
\def\H{{\mathcal H}}
\def\R{{\mathcal R}}

\def\B{{\mathcal B}}
\def\fH{\H^\fn}
\def\ann{{\operatorname{Ann}}}
\def\hom{{\operatorname{Hom}}}
\def\rad{{\operatorname{Rad}}}
\def\End{{\operatorname{End}}}
\def\dim{{\operatorname{dim}}}
\def\ind{{\operatorname{ind}}}
\def\res{{\operatorname{res}}}
\def\Rep{{\operatorname{Rep}}}
\def\soc{{\operatorname{soc}\:}}
\def\cosoc{{\operatorname{cosoc}\:}}
\def\ch{{\operatorname{ch}\:}}
\def\ba{{\underline{a}}}
\def\bb{{\underline{b}}}
\def\bc{{\underbar{c}}}

\def\ga{{\gamma}}
\def\Ga{{\Gamma}}
\def\la{{\lambda}}
\def\si{{\sigma}}
\def\De{{\Delta}}

\def\eps{{\epsilon}}

\begin{document}

\title{On representations of affine Hecke algebras of type $B$}

\author{Vanessa Miemietz}
\address{Universit\"at Stuttgart, Fachbereich Mathematik, Institut f\"ur
  Algebra und Zahlentheorie, Pfaffenwaldring 57, 70569 Stuttgart, Germany}
\email{miemieva@mathematik.uni-stuttgart.de}
\thanks{Research supported by the Studienstiftung des deutschen Volkes.}


\date{\today}



\begin{abstract}
Ariki's and Grojnowski's approach  to the representation theory of affine Hecke
algebras of type $A$ is applied to type $B$ with unequal parameters to obtain -- under
certain 
restrictions on the eigenvalues of the lattice operators -- analogous
multiplicity-one results and a classification of irreducibles with partial
branching rules as in type $A$. 
\end{abstract}
\maketitle

\section*{Introduction}

In this paper, the methods Ariki \cite{A:96} and Grojnowski \cite{G:99} developed for the representation theory
of affine Hecke algebras of type $A$ are applied to affine Hecke algebras of
type $B$.
The first section introduces the affine Hecke algebras $\H_n$  -- which are
the main objects of interest in this paper -- and their
subalgebras $\H_n^R$ which will be investigated in the last section. It is
explained how to use Clifford theory to exploit knowledge about one algebra
to obtain results about the other. 
The second through fourth sections closely follow Brundan and Kleshchev's paper \cite{BK:01} and 
informal lecture notes \cite{K:01} by Kleshchev that are now part of his book \cite{K:05}. The second section provides an affine version of the Mackey
Theorem and investigates the relation between induction and coinduction
functors. The third section introduces the concept of formal characters, which
are the main tool in understanding finite-dimensional irreducible modules for
the affine Hecke algebras of type $A$. The main results for $\H_n$ such as
irreducibility of the cosocle of certain induced modules and
multiplicity-freeness of the socle of certain restricted modules  are stated in the fourth
section. The fifth section contains some results in the cases where the methods used in type $A$ don't work.
After providing an overview of results on the affine Hecke algebra of
type $A$ in Section 6, we then give a one-to-one correspondence between
irreducibles in certain subcategories of the module category of $\H_n^R$
and irreducibles in the analogous subcategories of the module category of
$\H_n^A$ in the last section, yielding partial branching rules in those cases.

\section{The Algebras}

     Fixing an algebraically
closed field $F$ of characteristic not equal to two containing deformation
parameters $p$ and $q$ which are not roots of unity, we  define the affine
Hecke algebra of type $B_n$ for $n \geq 1$
to be the associative $F$-algebra  $\H_n$ on generators 
$$X_0^{\pm 1}, \dots, X_n^{\pm 1},  T_0, \dots, T_{n-1},$$ where the $T_i$
generate a finite Hecke algebra $\fH_n$ of type $B_n$  with relations
\begin{align*}\label{rels}
(1)&& (T_0-p)(T_0+p^{-1}) &= 0&\\
(2)&& (T_i-q)(T_i+q^{-1}) &= 0 &\hbox{ for } i \geq 1\\
(3)&& T_i T_{i+1} T_i &=
  T_{i+1} T_i T_{i+1}
  &\hbox{ for } i \geq 1 \\
(4)&&T_i T_j  &=  T_j T_i & \hbox{ for } |i-j|>1\\
(5)&&T_1 T_0 T_1 T_0 & =  T_0 T_1 T_0 T_1, & 
\end{align*}
and the $X_i^{\pm 1}$ generate a Laurent polynomial ring $\P_n$. Those two subalgebras are
subject to the mixed relations
\begin{align*}
(6)&&T_0 X_0 T_0 &= X_0X_1 & \\
(7)&&T_i X_j &= X_j T_i & \hbox{for } j \neq i,i+1\\ 
(8)&&T_i X_{i} T_i &= X_{i+1}& \hbox{for } i \geq 1.
\end{align*}

For $n=0$, we define $\H_0 :=F[X_0^{\pm1}].$ 

This is a deformation of the group algebra of the extended affine Weyl group
$W_n$ using the weight lattice of the general orthogonal group
$GO_{2n+1}(F)$, which is the subgroup of $GL_{2n+1}(F)$ respecting the
orthogonal form up to a scalar. $W_n$ is isomorphic to the
semidirect product of the finite Weyl group $W_n^\fn$ of type $B_n$ with generators
$s_0, s_1 \dots, s_{n-1}$ and relations
\begin{align*}
 s_i^2 &= 1 &\hbox{ for all } i\\
 s_i s_{i+1} s_i &=
  s_{i+1} s_i s_{i+1}
  &\hbox{ for } i \geq 1 \\
s_i s_j  &=  s_j s_i & \hbox{ for } |i-j|>1\\
s_1 s_0 s_1 s_0 & =  s_0 s_1 s_0 s_1, & 
\end{align*}
and the weight lattice of $GO_{2n+1}(F)$ which is the free abelian group on $n+1$
generators $X_0, \dots, X_n$, on which $W_n^\fn$ acts as in the mixed
relations for the affine Hecke algebra, substituting $s_i$ for $T_i$.
The actual affine Weyl group is the subgroup obtained from $W_n$ by omitting
the generator $X_0$ and adding the additional relation $s_0X_1s_0 =X_1^{-1}$
which in $W_n$ can be obtained from the first of the mixed relations.

The deformation of the affine Weyl group 
is then naturally a subalgebra of $\H_n$ generated by $X_1^{\pm 1}, \dots,
X_n^{\pm 1}$ and $T_0, \dots, T_{n-1}$ which we denote by $\H_n^R$. Here we
need an additional relation
$$X_1T_0= T_0X_1^{-1} + (p-p^{-1})(X_1+1)$$
which, in $\H_n$, can be derived from relation (6) since 
\begin{equation*}\begin{split}
X_1T_0 &= X_0^{-1}T_0X_0T_0T_0 \\
&=(p-p^{-1})X_0^{-1}T_0X_0T_0+X_0^{-1}T_0X_0\\
&=(p-p^{-1})X_1 +X_0^{-1}T_0^{-1}X_0 + (p-p^{-1})\\
&=(p-p^{-1})(X_1+1) + T_0X_0^{-1}X_1^{-1}X_0\\
&=T_0X_1^{-1} +(p-p^{-1})(X_1+1).
\end{split}\end{equation*}
The
commutative subalgebra generated by $X_1^{\pm 1}, \dots, X_n^{\pm 1}$ will be
denoted by $\R_n$. 

For a reduced expression $w=s_{i_1} \cdots s_{i_k}$ of an element $w \in W_n^\fn$, we define 
$T_w:=T_{i_1} \cdots T_{i_k}$. This does not depend on the choice of reduced expression and is therefore well-defined.

In \cite{L:89} Lusztig proves a general result on bases of affine Hecke algebras in the case where $p$ and $q$ are
distinct powers of the same deformation parameter $v_0$, but the proof doesn't
 rely on this and carries over to the general case, see
 \cite{Xi:94}. In our case this result gives the following two bases for
 $\H_n$:

\begin{equation}\label{basis1}
\left\{ X_0^{c_0}X_1^{c_1} \cdots X_n^{c_n}T_w 
\left| \substack{(c_0,\dots,c_n) \in \ZZ^{n+1},\\ w \in W_n^\fn} \right.\right\} \end{equation}

and

\begin{equation}\label{basis2}
\left\{ T_w X_0^{c_0}X_1^{c_1} \cdots X_n^{c_n}
\left| \substack{(c_0,\dots,c_n) \in \ZZ^{n+1},\\ w \in W^\fn_n} \right.\right\}. \end{equation}

All modules under consideration will be left modules that are finite-dimensional 
over $F$ and the category of such modules for an $F$-algebra $A$ will be denoted
by $A \mod$. For any affine Hecke algebra $\H$, Bernstein showed that its center $Z(\H)$ is exactly the set of
Laurent polynomials $f$ in its lattice that are invariant under the action of the
finite Weyl group on the lattice, see e.g.\ \cite{Xi:94}, \S 2.9. It is well-known that all irreducible representations of $\H_n$ are 
finite-dimensional, since $\H_n$ is finite-dimensional over its center, which by Dixmier's version of Schur's Lemma acts as a scalar on irreducible $\H_n$-modules.
 The Grothendieck group of the category $A \mod$ will be denoted by $K(A \mod)$ and for $M \in A \mod$ the corresponding element in $K(A \mod)$ will be written as $[M]$. For $M \in A \mod$, we denote the socle by $\soc M$ and the cosocle, i.e.\ the largest semisimple quotient, by $\cosoc M$.
If we have an automorphism $\psi$ of $A$, we will, for any $M \in
A \mod$ denote by $M^\psi$ the module obtained from $M$ by twisting the
action with $\psi$. This is equal to $M$ as an abelian group but the
operation of $A$ is now via the new multiplication $\diamond $ defined by \mbox{$a
\diamond m= \psi(a) m$} for $a \in A, m\in M$. The smallest integer $k$
such that  $M^{\psi^k} \cong M$ as an $A$-module will be called the order of
$\psi$ on $M$ whereas the order of $\psi$ (without specification of a module)
will denote the order of $\psi$ on $A$.

\subsection{Clifford theory}

We will use Clifford theory to move between modules for both algebras $\H_n^R$
and $\H_n$. This
idea to explore the interplay between different affine Hecke algebras
of the same isogeny class is originally due to Xi \cite{Xi:94} and has been worked
out in detail by Ram and Ramagge \cite{Ra:03}. In fact, Clifford theory works in a more general
setting, which has been studied by Dade in \cite{D:70}.  

\begin{lemma}\label{cliff}
Let $n$ be a natural number, $K$ an algebraically closed field of characteristic $p \geq 0$ with
$p \nmid n$. Let $A$ be a
$K$-algebra and let $B$ be a subalgebra of $A$, such that $A$ is free as a
$B$-module on basis $\{ x^s \mid 0 \leq s \leq n-1 \}$ for an
invertible element $x$ in $A$, and $\ZZ / n \ZZ$-graded,
i.e.\ $Bx^sBx^t = Bx^{s+t}$.
 Let $\psi: a \mapsto x^{-1}ax$ be conjugation with $x$, so 
$\psi (B) =B$. 
Let $M \in A \mod$ and let $N$ be an irreducible $B$-submodule of $\res^A_B M$.
Then the order $d$ of $\psi$ on $N$ divides $n$ and for $k:=n/d$  we have
$$\res^A_B M = \overset{d-1}{\underset{j=0}{\bigoplus}}N^{\psi^j}$$
and
$$\ind^A_B N = \overset{k-1}{\underset{j=0}{\bigoplus}}M_j$$
for irreducible and pairwise non-isomorphic modules $M_j$. Let $\si$ be an automorphism of $A$ with $\si \mid_B = \id_B$. All
$\si$-conjugates of $M$ occur as some $M_j$ in this decomposition, so
in particular, the
order of $\si$ on $M_j$ is less than or equal to $k$ for all $0 \leq j \leq k-1$.
\end{lemma}

\begin{proof}
Since $x^n \in B$, $N^{\psi^n} \cong N$ for any $B$-module $N$.
Now let $d$ be the smallest natural number such that $N^{\psi^d} \cong N$ and 
assume $d$ does not divide $n$. Then for $n=qd+r$, $N^{{\psi^{qd}}\psi^r} \cong
N$, i.e.\ $N^{\psi^r} \cong N$, but $r < d$, a contradiction. So, indeed $d$
does divide $n$.

Let $f: N \rightarrow N^{\psi^d}$ be an isomorphism and note that then $f^j: N
\rightarrow N^{\psi^{dj}}$ is also an isomorphism. In particular, since $x^n
\in B$, $f^k: N \rightarrow N^{\psi^{dk}}$ is a scalar multiple of
multiplication with $x^{-n}$, so by normalizing, we can assume that  $f^k$ is in
fact multiplication with $x^{-n}$.

Now take any irreducible $B$-submodule $N$ of $M$ and consider 
$\ind^{B'}_B N$, where we set $B':= \underset{0 \leq j \leq
  k-1}{\bigoplus}Bx^{jd}$.

\emph{Claim 1:} $\ind^{B'}_B N$ is a completely reducible $B'$-module, decomposing into a direct sum  of $k$ non-isomorphic $B'$-modules $L_i$, $i=0,\dots,k-1$, where each $L_i$ is isomorphic to $N$ as $B$-module.

\emph{Proof of Claim 1:}
As a $B$-module $\ind^{B'}_B N \cong \bigoplus_{j=0}^{k-1} x^{jd} \otimes
N.$
Let $\zeta$ be a primitive $k$-th root of unity and define $L_i$ to  be the subspace
of $\ind^{B'}_B N$ consisting of all elements
$a_i:=\sum_{j=0}^{k-1} \zeta^{ji} x^{jd} \otimes  f^j (a)$
where $a$ runs through $N$.
It is straightforward to check that $b \cdot a_i = (ba)_i$, so $L_i$ is a
$B$-submodule of $\ind^{B'}_B N$, giving a $B$-isomorphism between $L_i$ and $N$.
$L_i$ is also a $B'$-submodule of $\ind^{B'}_B N$. To see this, it suffices to
show $x^dL_i \subseteq L_i$ since $B'=
  \bigoplus_{j=0}^{k-1}B x^{jd}$  as a $B$-module. But for $a$ in $N$,
\begin{equation*}\begin{split}
x^d a_i 
&=\sum_{j=0}^{k-1} \zeta^{ji} x^{(j+1)d} \otimes f^j (a)\\
&=\zeta^{-i} \sum_{j=0}^{k-1} \zeta^{(j+1)i} x^{(j+1)d} \otimes f^{j+1}  f^{-1}(a)\\
&=\zeta^{-i} (f^{-1}(a))_i,
\end{split}\end{equation*}
 which is again an element of $L_i$.

Now let $0 \leq i,l \leq k-1$ and suppose $i \neq l$ but $L_i \cong L_l$, i.e.\ there exists a $B'$-module isomorphism 
$g :L_i \longrightarrow L_l$.  Observe that $ \res^{B'}_B L_i$  is isomorphic to
$\res^{B'}_B L_l$ via the isomorphism $\tilde g: \;a_i \mapsto a_l$ and this
is the only isomorphism up to a scalar by Schur's Lemma and the irreducibility
of $\res^{B'}_B L_i \cong N$. Since, if $g$ is an isomorphism, $\lambda g$ is,
for any $\lambda \neq 0 \in K$,
also an isomorphism, we can choose $g$ to coincide with the map $\tilde g$.
But then
\begin{equation*}\begin{split}
g(x^d a_i) &= \zeta^{-i} g((f^{-1}(a))_i)\\
&=\zeta^{-i}(f^{-1}(a))_l
\end{split}\end{equation*}
but also
\begin{equation*}\begin{split}
g(x^d a_i) &= x^d g(a_i)\\
&= x^d a_l\\
&= \zeta^{-l}(f^{-1}(a))_l
\end{split}\end{equation*}
whence we conclude that $i=l$, contrary to our assumption. This proves Claim
$1$. 

Now let $L$ be one of the irreducible $B'$-submodules of $\ind^{B'}_B N$ and
consider $\ind^A_{B'} L$. The set $\{ x^j \mid 0 \leq j \leq d-1 \}$ forms a
basis of $A$ as a $\ZZ / d \ZZ$-graded $B'$-module. The automorphism $\psi$ leaves $B'$
invariant as it leaves $B$ invariant and fixes $x^{jd}$ for $0 \leq j \leq
k-1$. Therefore we can twist any $B'$-module with $\psi$ and again obtain a
$B'$-module. Since $N^{\psi^j} \ncong
N$ for $j<d$, we also have $x^j\otimes L \cong L^{\psi^j} \ncong L$ for $j<d$.

\emph{Claim 2:} $\ind^A_{B'} L$ is an irreducible $A$-module.

 We have $$\ann_{B'} x^i \otimes L  \nsupseteq \underset{\substack{0 \leq j \le d-1 \\  j
   \neq i}}{\bigcap} \ann_{B'} x^j \otimes L  .$$ 
Otherwise an inclusion of the annihilators which are the kernels of the representations 
$$\rho_{\neq i}: \; B' \rightarrow \End_K(\underset{\substack{0 \leq j \le d-1\\  j \neq
  i}}{\bigoplus} x^j \otimes L)$$ and $$\rho_i: \; B' \rightarrow \End_K(x^i \otimes L)$$
would give rise to a projection on the side of the
images, which cannot happen since the $x^j\otimes L$ are pairwise
  non-isomorphic irreducible $B'$-modules. 

Let $0 \neq a = \sum_{j=0}^{d-1} x^j \otimes a_j$ with $a_j \in L$ be an
element of $\ind^A_{B'} L$.
Suppose $a_i \neq 0$. Then the left ideal $\ann_{B'}x^i \otimes a_i$ contains the maximal
 two-sided ideal $\ann_{B'} x^i \otimes L$. If it also contained $\underset{\substack{0 \leq j \le d-1\\  j
   \neq i}}{\bigcap} \ann_{B'} x^j \otimes L$, it would contain the two-sided
 ideal $\ann_{B'} x^i \otimes L + \underset{\substack{0 \leq j \le d-1\\  j
   \neq i}}{\bigcap} \ann_{B'} x^j \otimes L$, which, by the above, is strictly
 larger than $\ann_{B'} x^i \otimes L$ and thus equals $B'$ by the maximality
 of $\ann_{B'} x^i \otimes L$.

Therefore we can find a non-zero element such that $$ya =y x^i \otimes a_i = x^i \otimes
y' a_i \neq 0,$$ where $y'=x^{-i}yx^i$ and therefore we have an element $0 \neq b_i:= y' a_i \in L$ 
such that $x^i \otimes  b_i \in B'a$. Thus, $x^i \otimes L$ is a
$B'$-submodule of $B'a$ and $x^j \otimes L$ is contained in $Aa$ for all $j=0, \dots, d-1$. This proves Claim 2.

Now, what is left for us to prove is that the $M_{i}:= \ind_{B'}^A L_{i}$ are
also pairwise non-isomorphic and that all conjugates of $M$ by distinct powers
of $\si$ occur in the decomposition.

Frobenius reciprocity gives
 \begin{equation*}\begin{split}
 \hom_A(M_{i},M_{l}) &\cong
 \hom_{B'}(L_{i},\overset{d-1}{\underset{j=0}{\bigoplus}}x^j \otimes L_{l})
 \\&\cong \hom_{B'}(L_{i}, L_{l}) \oplus
 \overset{d-1}{\underset{j=1}{\bigoplus}} \hom_{B'}(L_{i}, x^j \otimes
 L_{l}).
\end{split}\end{equation*}
 $\hom_{B'}(L_{i}, L_{l})=0$ for $i \neq l$,
 according to Claim 1, while $$\overset{d-1}{\underset{j=1}{\bigoplus}}
 \hom_{B'}(L_{i}, x^j \otimes L_{l})\hookrightarrow
 \overset{d-1}{\underset{j=1}{\bigoplus}} \hom_{B}(\res^{B'}_B L_{i},
 \res^{B'}_B x^j \otimes L_{l}) $$ but the latter is just
 $\overset{d-1}{\underset{j=1}{\bigoplus}} \hom_{B}(N, N^{\psi^j})$ which is zero by hypothesis.
So the $M_{i}$ are pairwise non-isomorphic.

Now we have $\ind_B^A N \cong \overset{k-1}{\underset{i=0}{\bigoplus}}M_{i}$
and $\res^A_B M_{i}\cong \overset{d-1}{\underset{j=0}{\bigoplus}} x^j \otimes
N$. In particular, this implies the asserted decomposition of $M$, since $M$
is isomorphic to one of the $M_i$ by Frobenius reciprocity. 
As $\si$ fixes $B$ pointwise, $S^{\si^l} = S$ for all $B$-modules $S$ and all $l \in \ZZ$. So
 \begin{equation*}\begin{split}
K &\cong \hom_B(N, \res^A_B M_{i}) \\
&\cong \hom_B(N, \res^A_B M_{i}^{\si^l})\\ 
&\cong \hom_A(\overset{k-1}{\underset{j=0}{\bigoplus}}M_{j}, M_{i}^{\si^l})
\end{split}\end{equation*} 
gives an inclusion of $\{ M_{i}^{\si^l} \mid l \in \ZZ \}$
into $\{M_{j} \mid 0\le j \le k-1  \}$ for a fixed $i$. Therefore the order of $\si$ on $M_{i}$ is at most $k$ for any $i$.
\end{proof}

Now, in order to apply this to our situation,
 observe that $X_0^2 \prod_{i=1}^n X_i$ is central in $\H_n$. To verify
this, note that every generator $T_i$ commutes with all but $X_i$ and
$X_{i+1}$. But, for $i \geq 1$, $T_iX_iX_{i+1}= X_{i+1}T_i^{-1}X_{i+1} =
X_iX_{i+1}T_i.$
Similarly, $X_0^2X_1T_0 = X_0X_1T_0^{-1}X_0X_1 =T_0X_0^2X_1, $ 
 so $X_0^2
\prod_{i=1}^n X_i$ is indeed central and invertible 
 and therefore
acts as a nonzero scalar $\mu_M$ on an irreducible module \mbox{$M \in \H_n \mod$.} Thus,
 the action of $\H_n$ on an irreducible module $M \in \H_n \mod$
factors over the quotient algebra \mbox{$\H_n^{\mu_M} := \H_n/(X_0^2 \prod_{i=1}^n
X_i-\mu_M)$.}

\begin{lemma}
$\H_n^{\mu_M}$ is a $\ZZ/2 \ZZ$-graded $\H_n^R$-module with basis $\{1,X_0\}$.
\end{lemma}

\begin{proof}
We first show that $\H_n^R$ can be identified with a subalgebra of $\H_n^{\mu_M}$.
Since $X_0^2 \prod_{i=1}^nX_i$ is central the ideal generated by $X_0^2 \prod_{i=1}^n
X_i-\mu_M$ is already contained in and thus equal to the right ideal $(X_0^2 \prod_{i=1}^n
X_i-\mu_M)\H_n$ which is generated over $F$ by 
\begin{equation*}\begin{split}
(X_0^2 &\prod_{i=1}^n
X_i-\mu_M) X_0^{c_0}X_1^{c_1} \cdots X_n^{c_n}T_w\\&=
 (X_0^{c_0+2}X_1^{c_1+1} \cdots X_n^{c_n+1}-\mu_M X_0^{c_0}X_1^{c_1} \cdots X_n^{c_n})T_w
\end{split}\end{equation*}
for $(c_0,\dots,c_n) \in \ZZ^{n+1}, $ and $w
\in W_n^\fn $.

No finite linear combination of those can have degree $0$ in $X_0$, therefore
the intersection of the ideal generated by $(X_0^2 \prod_{i=1}^n
X_i-\mu_M)$ with $\H_n^R$ is zero and we can view $\H_n^R$ as a subalgebra of 
$\H_n^{\mu_M}$.

As $X_0$ commutes with all generators $X_j$ and all generators $T_j$ except
$T_0$ and 
\begin{equation*}\begin{split}
T_0X_0 &= T_0^{-1}X_0 + (p-p^{-1})X_0\\
&=X_0T_0X_1^{-1} + (p-p^{-1})X_0\\
&=X_0(T_0X_1^{-1} + (p-p^{-1}))
\end{split}\end{equation*}
we see that $\H_n^R X_0 \subseteq X_0 \H_n^R$, which by the same argument as
for $\H_n^R$ above has no nontrivial intersection with the ideal generated  by
$(X_0^2 \prod_{i=1}^n
X_i-\mu_M)$ and can therefore be viewed as contained in $\H_n^{\mu_M}$.
Certainly $\H_n^R \cap \H_n^R X_0 = \{0\}$, so we have an
$\H_n^R$-submodule of $\H_n^{\mu_M}$ which is isomorphic to $\H_n^R \oplus
\H_n^R X_0 $.
But \mbox{$X_0 \H_n^R X_0 \subseteq \H_n^R +(X_0^2 \prod_{i=1}^n
X_i-\mu_M) \H_n $,} as we see by considering 
that $X_0^2 = \mu_M(\prod_{i=1}^n
X_i^{-1}) \in \H_n^{\mu_M}$ and $X_0T_0X_0 = T_0^{-1}X_0^2X_1$.
Thus $\H_n^{\mu_M} \cong \H_n^R \oplus
\H_n^R X_0 $ as a left $\ZZ/2 \ZZ$-graded $\H_n^R$-module.
\end{proof}

Set 
$$\psi:\  \H_n \rightarrow \H_n: \ \ h \mapsto X_0^{- 1} h
X_0$$ and $$\begin{array}{crll} \si: \ \H_n \rightarrow \H_n: \quad & T_i& \mapsto
T_i&\\ & X_0& \mapsto -X_0& \\ &X_i& \mapsto X_i &\textrm{ for $i \geq 1$.} \end{array}$$
Both automorphisms leave the ideal generated by $(X_0^2 \prod_{i=1}^n
X_i-\mu_M)$ invariant.
For $\psi$, this follows from the commutativity of $\P_n$ and for
$\si$ from the fact that $(-X_0)^2=X_0^2$. 
Thus $\psi$ and $\si$ define automorphisms on $\H_n^{\mu_M}$.
Then $\psi$ leaves $\H_n^R$ invariant since $$X_0^{-1}T_0X_0 = X_1T_0^{-1} $$
and $\si$ fixes $\H_n^R$
pointwise since it is the identity on the generators of $\H_n^R$. Applying Lemma \ref{cliff}, the restriction to $\H_n^R$ of an irreducible
$\H_n^{\mu_M}$-module $M$ splits only if $M^{\si} \cong M$. 
Recall that $\operatorname{char} F \neq 2$.

\begin{lemma}\label{cliff2}
 $M^{\si} \cong M$ only if $-1$ occurs as an eigenvalue for some
  $X_j$, \mbox{$j=1, \dots, n$.} 
\end{lemma}

\begin{proof} 
The element $X_0 \underset{1 \leq i \leq n}{\prod} (1+X_i)$ is central in $\H_n$ as 
\begin{equation*}\begin{split}
X_0(1+X_1)T_0 &= X_0T_0 + X_0X_1T_0\\
&= T_0^{-1}X_0X_1 + T_0X_0 + (p-p^{-1})X_0X_1\\
&=T_0(X_0+X_0X_1)
\end{split}\end{equation*}
and 
\begin{equation*}\begin{split}
(X_i+1)(X_{i+1}+1)T_i&= X_iX_{i+1}T_i + X_iT_i + X_{i+1}T_i +T_i\\
&=T_iX_iX_{i+1} +T_i^{-1}X_{i+1} + T_iX_i \\& \quad +(q-q^{-1})X_{i+1} +T_i\\
&=T_i(X_i+1)(X_{i+1}+1)
\end{split}\end{equation*}
and all other factors commute with $T_i$ anyway.
 If $-1$ does not occur as an eigenvalue for any of the $X_i$, $1 \leq i \leq n$, 
$X_0 \underset{1 \leq i \leq n}{\prod} (1+X_i)$ acts by a nonzero scalar on $M$ and 
by its negative on $M^{\si}$, so the two are not isomorphic.
\end{proof}

If $-1$ occurs as an eigenvalue of some $X_i$ on $M$, it can indeed happen that $M^{\si} \cong M$ but this will not always be the case.
\begin{example} Consider the 2-dimensional module $M$ for $\H_1$ on which the generators $T_0, X_0,X_1$ act by the matrices
$$ \begin{array}{ccc}
 \begin{pmatrix} p & 0\\ 0& -p^{-1}
\end{pmatrix}  , &
\begin{pmatrix} 0& a_0\\ a_0 & 0
\end{pmatrix} ,&
\begin{pmatrix} -1 & 0\\ 0 & -1
\end{pmatrix} 
\end{array}   $$
respectively. This is obviously irreducible but splits upon restriction to $\H_1^R$, the isomorphism 
between $M$ and $M^\si$ being given by multiplication with the matrix $\begin{pmatrix} 1&0\\0&-1 \end{pmatrix}$.
On the other hand, there is a two-dimensional representation for $\H_2$
where $T_0$ and $T_1$ act as
$$ \begin{array}{ccc}
 \begin{pmatrix} \frac{(p-p^{-1})q^2}{(q^2+1)} &
\frac{p^4q^2+q^4p^2+p^2+q^2}{pq^2(p^2-1)(1+q^2)}
\\ \frac{q^2(p-p^{-1})}{(q^2+1)}&  \frac{(p-p^{-1})}{(q^2+1)}
\end{pmatrix}  & \hbox{and} &
 \begin{pmatrix} -q^{-1} & 0\\ 0& q
\end{pmatrix}  
\end{array}$$
and $X_0,X_1$ and $X_2$ act as
$$ \begin{array}{cccc}
\begin{pmatrix} a_0& 0\\ 0&-a_0q^2
\end{pmatrix} ,&
\begin{pmatrix} -q^2 & 0\\ 0 & -q^{-2}
\end{pmatrix}& \hbox{and}&
\begin{pmatrix} -1 & 0\\ 0 & -1
\end{pmatrix} 
\end{array}   $$ respectively,
which is also irreducible and which remains irreducible when restricted to $\H_2^R$.
\end{example}
Since our main goal is to understand finite dimensional irreducible modules for affine Hecke algebras
by looking at the action of the lattice, we will generally work with the
algebra $\H_n$ as opposed to the algebra $\H_n^R$. We do not know whether the
action of the lattice uniquely determines irreducibles for $\H_n$ but at
least it does in all examples we have been able to compute, contrary to
$\H_n^R$ where the above is an immediate counterexample.

\subsection{Notation and computations}
The affine Hecke algebra of type $\tilde A_n$ using the weight lattice of $GL_n$ is
naturally embedded in $\H_n$ as the subalgebra generated by 
$T_1, \dots, T_{n-1}$ and $X_1^{\pm 1}, \dots,X_n^{\pm 1}$ which will be denoted by 
$\H_n^A$. A lot of work has been done on this algebra, upon which we will
heavily rely in the following. We will review the theory in type $A$ when
needed and point out similarities and differences along the way.

To simplify notation, we now introduce some abbreviations. 
A simultaneous 
eigenvector for a set of lattice operators $X_{i_1},\dots, X_{i_r}$ with
respective eigenvalues $a_{i_1},\dots, a_{i_r}$ will be called an
$(a_{i_1},\dots, a_{i_r})$-eigenvector for $X_{i_1},\dots,
X_{i_r}$. Analogously, we will refer to the subspace of all 
$(a_{i_1},\dots, a_{i_r})$-eigenvectors for $X_{i_1},\dots,X_{i_r}$ as the
$(a_{i_1},\dots, a_{i_r})$-eigenspace for $X_{i_1},\dots,X_{i_r}$. If the 
specification ``for
 $X_{i_1},\dots,X_{i_r}$'' is omitted, the tuple $(a_{i_1},\dots, a_{i_r})$ 
necessarily 
 consists of as many entries as there are lattice operators and the lattice 
operators are 
 taken in order, i.e.\ the $(a_0,a_{1},\dots, a_{n})$-eigenspace in a module $M \in
 \H_n \mod$ denotes 
 the $(a_0,a_{1},\dots, a_{n})$-eigenspace for $X_0,X_1,\dots,X_n$ .
The same conventions will be used for generalized eigenvectors and eigenspaces.

We will write $T_kT_{k+1} \cdots T_l$ for $k \leq l$ or $T_kT_{k-1} \cdots T_l$ for 
$k \geq l$ and $T_k \cdots T_0 \cdots T_l$ as $T_{k,l}$ and $T_{k,0,l}$ respectively, and adopt the analogous convention for $s_{k,l}$ and $s_{k,0,l}$.

Next, we state some technical lemmas which are easily checked by direct
computations using the defining relations.

\begin{lemma}\label{intertwining}
Let $N \in \H_n \mod$.
\begin{enumerate}
\item Let $j \geq 1$ and let $u \in N$ be an $(a,b)$-eigenvector for $X_j,X_{j+1}$ where $a,b \in
  F$ satisfy $a^{-1}b \notin \{ q^{\pm 2}\}$. Then $(T_j- a^{-1}bT_j^{-1})u$
  is a $(b,a)$-eigenvector for $X_j,X_{j+1}$.
\item Let $u \in N$ be an $(a_0,a)$-eigenvector for $X_0,X_1$, where $a \notin
  \{ p^2, 1  \}$. Then $(T_0 - a T_0^{-1})u$ is an $(a_0a,a^{-1})$-eigenvector
  for $X_0, X_1$.
\end{enumerate}
\end{lemma}

The restrictions on $a$ and $b$ guarantee that the vectors $(T_j- a^{-1}bT_j^{-1})u$ resp.\ $(T_0 - a T_0^{-1})u$ are nonzero.

We will now consider the behavior of  elements in
$\H_n$ when we move lattice elements from one side to the other but we first need
to define the Bruhat order on $W_n^\fn$: For $x,y \in W_n^\fn$, $x
< y$ if and only if there exists a reduced expression $y=u_1 \cdots u_k$
where $u_j \in \{s_i \mid i\in I\}$ for $1 \leq j \leq k$ and a subsequence \mbox{$1
\leq m_1 < \cdots < m_l \leq k$} such that $x= u_{m_1} \cdots u_{m_l}$. This
defines a partial order on $W_n^\fn$ which is compatible with the length
function. The following lemma is also straightforward.

\begin{lemma}\label{righttoleft}
Let $w \in W_n^\fn$ and $0 \leq i \leq n$. Then 
$$T_w X_i \in X_{w(i)} T_w + \sum_{\tilde{w} <
  w}\P_nT_{\tilde{w}}$$
where $X_{w^{-1}(i)}$ is the element $w^{-1}X_iw$ when we look at the action of the finite Weyl group on the lattice. 
\end{lemma}

\section{Mackey Filtration and Duality}

In a finite Coxeter group, a parabolic subgroup is generally defined as
being conjugate to a standard parabolic subgroup, which is the
subgroup generated by a subset of the Coxeter generators $s_i$. For $W_n^\fn$
this means we take a subset 
\begin{equation*}\begin{split}
I &=\{i_1,i_1+1,\dots,i_1+{r_1}-1,i_2,
\dots,i_2+r_2-1, \dots, i_l, \dots, i_l+r_l-1 \\& \quad \quad \quad \quad \quad
\quad \quad \quad \quad \quad \quad \quad \mid i_k +r_k< i_{k+1} \quad \forall 1 \leq
k \leq l-1 \}\\
& \subseteq \{0,1, \dots n-1\}
\end{split}\end{equation*}
 and obtain $W_I^\fn:=
\langle s_i \mid i \in I \rangle$ which is isomorphic to $W_{r_1}^\fn \times
S_{r_2} \times \cdots \times S_{r_l}$ if $i_1=0$ and to $S_{r_1} \times
S_{r_2} \times \cdots \times S_{r_l}$ if $0 \notin I$.

We generalize this concept to (extended) affine Weyl groups not by taking a subset of the
Coxeter generators which would yield finite Coxeter groups of infinite index
in the affine Weyl group --- making induction difficult to handle --- but by taking the
semidirect product of a (standard) parabolic subgroup of the finite Weyl group with the
full translation lattice. Thus we define $W_I := \langle s_i,X_j \mid i \in I, j
\in \{0, \dots, n\} \rangle$.

Analogously, the parabolic subalgebra $\H_I^\fn$ of $\H_n^\fn$ is defined as 
the subalgebra generated by $\{  T_i \mid i \in I \}$ and the parabolic subalgebra
$\H_I$ of $\H_n$ is the subalgebra generated by $\{  T_i, X_j^{\pm 1} \mid i \in
I , j \in \{0, \dots, n\} \}$. 
We will also write $(m_0, \dots, m_l)$ with $m_0 \geq 0$ and $m_i \geq 1$ for $i \geq 1$ where $m_0+ \cdots +m_l=n$ for
$$I=\{0,1,\dots, m_0-1,m_0+1, \dots, m_0+m_1-1, \dots, \sum_{i=0}^{l-1} m_i
+1,\dots,n-1\}$$ and write $\H_{m_0, \dots, m_l}$ for $\H_I$ which is isomorphic
to \mbox{$\H_{m_0} \otimes \H^A_{m_1}\otimes \cdots \otimes\H^A_{m_l}$.}
In case $m_0=0$, $\H_{0, m_1, \dots, m_l} \cong F[X_0^{\pm1}] \otimes \H_{m_1, \dots, m_l}^A \cong F[X_0^{\pm1}] \otimes \H^A_{m_1} \otimes \H^A_{m_2}\otimes \cdots \otimes\H^A_{m_l}$ denotes the  tensor product of $F[X_0^{\pm1}]$ with the parabolic subalgebra of type $A$ corresponding to \begin{equation*}\begin{split}I&=\{1,\dots, m_1-1,m_1+1, \dots, m_1+m_2-1, \dots, \sum_{i=1}^{l-1} m_i +1,\dots,n-1\}.\end{split}\end{equation*}


We will generally abbreviate induction and restriction functors between
parabolic subalgebras as $\ind^I_J := \ind^{\H_I}_{\H_J}$, $\res^I_J :=
\res^{\H_I}_{\H_J}$, \mbox{$\ind^{n_0, \dots, n_k}_{m_0, \dots, m_l} := \ind^{\H_{n_0, \dots, n_k}}_{\H_{m_0, \dots, m_l}}$} and $\res^{n_0, \dots, n_k}_{m_0, \dots, m_l} :=
\res^{\H_{n_0, \dots, n_k}}_{\H_{m_0, \dots, m_l}}$ for $J \subseteq I$ or
$({m_0, \dots, m_l}) \subseteq (n_0, \dots, n_k)$ respectively. If we induce directly from a parabolic subalgebra of type $A$, we will always use the full expression.

For a parabolic subgroup $W_I^\fn$ of $W_n^\fn$, there are
distinguished left  and right coset representatives of minimal length, the sets
of which will be denoted by $D_I$ and $D_I^{-1}$ respectively. For parabolic
subgroups $W_I^\fn$ and $W_J^\fn$ of $W_n^\fn$, \mbox{$D_{I,J}:= D_I^{-1} \cap D_J$}
is then the set of distinguished minimal length $(W_I^\fn,W_J^\fn)$-double coset
representatives. An account of this, including the following three properties
of distinguished double coset representatives, can be found in \cite{GP:00},
Chapter 2.1.

\begin{enumerate}
\item For $x \in D_{I,J}$,
$W^\fn_I \cap x W^\fn_J x^{-1} =: W^\fn_{I \cap {x} J}$ and
$x^{-1} W^\fn_I x \cap W^\fn_J =:W^\fn_{{{x^{-1}}}I \cap J}$ are parabolic subgroups of $W^\fn_n$. This defines subsets $I \cap {x} J$ and ${{x^{-1}}}I \cap J$ of $\{0, 1 , \dots, n-1\}$.

\item For $x \in D_{I,J}$, the map 
\begin{equation*}\begin{split}
W^\fn_{I\cap {x}J}& \rightarrow W^\fn_{{{x^{-1}}}I\cap J}\\
w & \mapsto x^{-1}wx
\end{split}\end{equation*}
defines a length preserving isomorphism. 

\item For $x \in D_{I,J}$, every
$w \in W^\fn_I x W^\fn_J$ can be written as $w = u x v$ for unique
elements $u \in W^\fn_I$ and $v \in W^\fn_J \cap D^{-1}_{x^{-1}I\cap J}$.
Moreover, $W^\fn_J \cap D^{-1}_{x^{-1}I\cap J}$ is the set of minimal 
length
right coset representatives of $W^\fn_{x^{-1}I\cap J}$ in $W^\fn_J$. 
\end{enumerate}

\begin{lemma}\label{ugly}
For $x \in D_{I,J}$, the subspace $\fH_{I}T_x\fH_{J}$ has basis $\{
T_w \mid w \in W^\fn_{I} x W^\fn_{J}\}.$

\end{lemma}
\begin{proof}
The proof of this is based on a standard argument using coset representatives and will be omitted.
\end{proof}

\begin{lemma}\label{BIJx}
For $x \in D_{I,J}$ the subspace
$\H_I T_x\fH_{J}$ of
$\H_n$ has basis $$B_{I,J}^x:=\{ X^{c_0}_0 \cdots X^{c_n}_n T_w \mid (c_0, \dots,c_n)\in \mathbb{Z}^{n+1},
w \in W^\fn_{I}x W^\fn_{J}  \} .$$
Moreover, as a vector space, $$\H_n = \underset{x \in D_{I,J}}{\bigoplus}\H_{I}T_x\fH_{J}. $$
\end{lemma}

\begin{proof} Analogous to \cite{BK:01}, Lemma 2.5.
\end{proof}

 We now fix some total order $\prec$ refining the Bruhat order $<$ on $D_{I, J}$. For $x \in D_{I, J}$, set
\begin{equation}\begin{split}\label{filt}
\B_{\preceq x} &= \bigoplus_{y \in D_{I, J},\ y \preceq x} 
\H_I T_y \fH_J,\\
\B_{\prec x} &= 
\bigoplus_{y \in D_{I, J},\ y \prec x} \H_I T_y \fH_J,\\
\B_{x} &= \B_{\preceq x} / \B_{\prec x}.
\end{split}\end{equation}
This defines a filtration of $\H_n$ as an $(\H_I,\H_J)$-bimodule, since it
follows from Lemma \ref{righttoleft} that we can move lattice elements from
the right to the left and only create terms that are smaller in the Bruhat
order and therefore lower in the filtration. 
Property (ii) of double coset representatives above implies  that for each $x \in D_{I, J}$, 
there exists an algebra isomorphism
\begin{equation*}\begin{split}
\phi_{x^{-1}}:\H_{I\cap x J} \rightarrow \H_{{x^{-1}}I\cap J}: T_w & \mapsto T_{x^{-1}wx}\\
X_j & \mapsto X_{x^{-1}(j)}
\end{split}\end{equation*}
 for $0 \leq j\leq n-1$ and $w \in W^\fn_{I \cap xJ}$. Note that $X_{x^{-1}(j)}$ is not necessarily of the form $X_i$ but can be the product of several polynomial generators.
 For $N \in \H_{{x^{-1}}I\cap J} \mod$, 
 ${^x N}$ will denote the $\H_{I\cap x J}$-module obtained by pulling back the
 action through $\phi_{x^{-1}}$.

Now we can prove an affine version of the Mackey theorem, which differs from
``classical'' Mackey theorems in that it does not give a direct decomposition
but only a filtration.

\begin{theorem} {\rm (``Mackey Theorem'')}
\label{TMackey}
Let $M \in \H_J \mod$. Then
$\res^{n}_{I}\ind^{n}_{ J} M$ admits a filtration with subquotients 
isomorphic to
$\ind_{{I\cap x J}}^{I}{^x}
(\res_{{{x^{-1}}I\cap J}}^{ J}M),$
one for each $x\in D_{I, J}$. 
The subquotients can be taken in any order refining the 
Bruhat order on $D_{I, J}$,
in particular, since the double coset representative  $1$ is the smallest
element in the ordering, $\ind_{{I\cap  J}}^{I}
\res_{I\cap J}^{ J}M$ appears as a submodule.
\end{theorem}

\begin{proof} 
We already have a filtration of $\H_n$ as an $(\H_I,\H_J)$-bimodule given in
(\ref{filt}). Thus $\res^{n}_{I}\ind^{n}_{ J} M ={_{\H_I}}\H_n \otimes_{\H_J}
M$ inherits a filtration with subquotients isomorphic to $\B_x \otimes_{\H_J}
M$, the $x \in D_{I,J}$ taken in any order refining the Bruhat order on $D_{I,J}$. 
Now 
\begin{equation*}\begin{split}
\ind_{{I\cap x J}}^{I}{^x}
(\res_{{{x^{-1}}I\cap J}}^{ J}M) &= \H_I\otimes_{\H_{I\cap x J}}
{^x}({}_{\H_{{x^{-1}}I\cap J}}\H_J \otimes_{\H_J} M)\\
& \cong \H_I\otimes_{\H_{I\cap x J}} {^x  \H_J} \otimes_{\H_J} M,
\end{split}\end{equation*}
thus it suffices to show that 
$\H_I\otimes_{\H_{I\cap x J}} {^x  \H_J} \cong \B_x.$
In order to do this, 
we define a bilinear map 
\begin{equation*}\begin{split}
\H_I \times  {}^x \H_J& \rightarrow \B_x\\
(h,h')& \mapsto h T_x h' +\B_{\prec x}.
\end{split}\end{equation*} 
Since, for $w \in \H_{I\cap x J}$, $ T_w T_x = T_{wx} =  T_{xx^{-1}wx} =
T_x T_{x^{-1}wx} $, and for all $j$, we have $X_jT_x=T_{x^{-1}}^{-1}T_{x^{-1}}X_jT_x = T_{x^{-1}}^{-1}X_{x^{-1}(j)} \in T_xX_{x^{-1}(j)} + \B_{\prec x}, $ this map is $\H_{I\cap x J}$-balanced and therefore induces
a map $\H_I \otimes_{\H_{I\cap x J}} {}^x \H_J \rightarrow \B_x.$

By Lemma \ref{ugly} and  property (iii) above it,
 a basis of $\H_I \otimes_{\H_{I\cap x J}} {}^x \H_J$ is given by  
$$ \{ X^{c_0}_0 \cdots X^{c_n}_n T_u \otimes T_v \mid (c_0, \dots,c_n) \in \mathbb{Z}^{n+1},
u \in W^\fn_{I},v \in  W^\fn_{J}\cap D^{-1}_{x^{-1}I\cap J}  \}, $$
the elements of which map to a basis of $\B_x$ by Lemma \ref{BIJx}, whence the
map is actually an isomorphism.\end{proof}

In general, for $N \in \H_n \mod$ and $M \in \H_I \mod$, by Frobenius reciprocity
\begin{equation}\label{frob}\hom_{\H_n}(\H_n \otimes_{\H_I} M,N) \cong \hom_{\H_I}(M, \res^{\H_n}_{\H_I}N)\end{equation} and 
\begin{equation}\label{otherfrob}\hom_{\H_I}(\res^{\H_n}_{\H_I} N,M) \cong \hom_{\H_n}(N,\hom_{\H_I}(\H_n,M)).\end{equation}
We would like to express the coinduced module $\hom_{\H_I}(\H_n,M)$ in terms of an induced module.
Hence, for the rest of this section, we fix a subset $I$ of \mbox{$\{0,1, \dots, n-1\}$.}
Let $d$ be the longest element of $D_{I,I}$.

\begin{lemma}  
Let $I=(m_1, \dots, m_l)$.
 Then the longest double coset representative $d$ in $D_{I,I}$ is an involution and $I \cap dI=I\cap d^{-1}I=I$.
\end{lemma}

\begin{proof}
Let $w_0$ be the longest element of $W_n^\fn$. This element has to map the
short root $X_1$ in the basis 
of the root system to a short root in the inverse of this basis which is
$X_1^{-1}, X_2^{-1}X_1, \dots, X_n^{-1}X_{n-1}$. The only short root there is
$X_1^{-1}$ and since $w_0$ acts as an isometry it follows that it sends each root
to its negative, hence also each $X_i$ to $X_i^{-1}$.
Set $k_i= \sum_{j=1}^im_j$ and let $w_{0,I}$ be the longest element in $$W_I^\fn \cong
W_{m_1} \times S_{m_2} \times \cdots \times S_{m_l},$$ which is
the element sending $X_1, \dots, X_{m_1}$ to their inverses as above and
reversing the 
orders of $X_{k_i+1}, \dots, X_{k_i+m_{i+1}}$ for $1
\leq i \leq l-1$, as the element reversing the order of the numbers $1, \dots, m_i$ is the
longest element in the symmetric group $S_{m_i}$.

For the longest distinguished left coset representative $\tilde d$ in $D_I$, we have 
\mbox{$w_0 = \tilde d w_{0,I}$} by the additivity of lengths for distinguished coset
representatives. Since $w_{0,I}$ is equal to its inverse, we obtain $\tilde d =
w_0 w_{0,I}$.
Computing the action of $\tilde d$ on the $X_i$, we see that it leaves $X_1,
\dots, X_{m_1}$ invariant and maps the ordered sets $(X_{k_i+1}, \dots,
X_{k_i+m_{i+1}})$ to $ (X_{k_i+m_{i+1}}^{-1},
\dots,X_{k_i+1}^{-1})$ for $1 \leq i \leq l-1$.
From this presentation it is easy to see that ${\tilde d}^{-1}$ is equal to
$\tilde d$, whence it is also the longest distinguished right coset
representative and therefore the longest element $d$ in $D_{I,I}$.
By direct computation it follows that $ds_id=s_i$ for $0 \leq  i \leq  m_1-1$
and $d s_{k_i+j}  d = s_{k_{i+1}+1-j}$ for $1 \leq j
\leq m_{i+1}$ and $1 \leq i \leq l-1$, which shows that
 $I \cap dI=I\cap d^{-1} I=I$.
\end{proof}

By property (ii) of distinguished double coset representatives, there is an isomorphism
\begin{equation}\label{meme}
\phi = \phi_{d^{-1}}:\H_I \rightarrow \H_I,
\end{equation}
 and for $M \in \H_I \mod$, we denote by
${^d M}$ the $\H_I$-module obtained by twisting the action with $\phi$.

We will need a homomorphism $\theta:\H_n \rightarrow {^d \H_I} $ of 
$(\H_I,\H_I)$-bimodules where the right action of $\H_I$ on ${^d \H_I}$ is the usual (untwisted) one. This homomorphism is given by first projecting 
$\H_n \longrightarrow \B_d$ in (\ref{filt}) and then applying the isomorphism
 of $\B_d \rightarrow \H_I \otimes_{\H_{I}} {^d} \H_I \cong {^d} \H_I $ 
 given in the proof of Theorem \ref{TMackey}. Explicitly, this homomorphism is given by 
$$\theta(X T_w) = \left\{
\begin{array}{ll}
\phi(X) T_{d^{-1}w} &\hbox{if $w \in d W_I^\fn$,}\\
0&\hbox{otherwise,}
\end{array}\right.$$
for $X \in \P_n,  w \in W_n^\fn$.
Then the following holds.

\begin{lemma}\label{joe3}

The map 
\begin{equation*}\begin{split}
f:\H_n &\rightarrow \hom_{\H_I}(\H_n,{{}^d \H_I}) \\
h &\mapsto (h \theta : t \mapsto \theta(th))
\end{split}\end{equation*}
 is an isomorphism of $(\H_n,\H_I)$-bimodules.
\end{lemma}

\begin{proof}
First, we need to show that $f$ is an $(\H_n,\H_I)$-bimodule homomorphism. 
So, we check $f(h)(t)=\theta(th) = h\theta (t) =hf(1)(t)$ for $h \in \H_n$ and $f(h')(t) =\theta(th') = \theta(t)h'= f(1)h'(t)$ for $h' \in \H_I$. 

Since, as left $\H_I$-module, $^d \H_I$ is isomorphic to $\H_I$ and $\H_n$ is a free left 
$\H_I$-module on basis $\{T_w \mid w \in D_I^{-1}\}$, 
the set $K:=\{\psi_w \mid w \in D_I^{-1} \}$, where 
\mbox{$\psi_w: \H_n \longrightarrow ^d \H_I: \:\: \psi_w(T_u)=\delta_{u,w}$} for $u \in D_I^{-1}$,
 is a basis for $\hom_{\H_I}(\H_n,{^d \H_I})$ as a free right $\H_I$-module.
 
So, we'll be done if we can give a basis for $\H_n$ as free right $\H_I$-module that is mapped to $K$.
But since $\H_I=\H_I^\fn \P_n$, a basis for $\H_n^\fn$ as free right $\H_I^\fn$-module is automatically a basis for $\H_n$ as free right $\H_I$-module.
Therefore, we study the restrictions $\theta ' :=\theta |_{\H_n^\fn}$ and $f' :=f |_{\H_n^\fn}$ to $\H_n^\fn$ of the above homomorphisms and want to construct a basis for $\H_n^\fn$ as free right $\H_I^\fn$-module such that the basis elements map to the $\psi_w ' :=\psi_w |_{\H_n^\fn}$. 
A basis like this can be found if we can show that $f'$ is an isomorphism of $(\fH_n, \fH_I)$-bimodules.
Suppose $f'$ is not an isomorphism, then there is a nonzero $h$ in its kernel,
i.e.\ $(f'(h))(t) = \theta'(th) = 0$ for all $t \in \fH_n$. Now write $h =
\sum_{y \in D_I, l(y) \leq l(x)} T_y h_y$ for some $x \in D_I, h_y \in \H_I$. If $x=d$, we know that $f'(h)(1)=\theta(h)=h_d \neq 0$, so we can use downward induction on $l(x)$ to show $h=0$. 
If $l(x) \le l(d)$ we can find a transposition $s$ such that $sx \in D_I$ and $l(sx) \geq l(x)$, so
$T_s h= \sum_{y \in D_I, l(y) \leq l(sx)} T_y h'_y $
and $h'_{sx} =h_x \neq 0$, so by the inductive assumption $\theta'(\fH_nh) =
\theta'(\fH_nT_sh) \neq 0$, whence $f'$ is indeed an isomorphism and we're
done. \end{proof}

\begin{mycor}\label{c1}
For $M \in \H_I \mod$, there is a natural isomorphism
$$\hom_{\H_I}(\H_n,{ M})
\cong
\H_n \otimes_{\H_I} {^dM}$$ of $\H_n$-modules.
\end{mycor}

\begin{proof}
Analogous to \cite{BK:01}, Corollary 2.12.
\end{proof}

On $\H_n$, 
we can define an 
antiautomorphism $\tau$ defined on the generators as follows:
\begin{equation*}\begin{split}
\label{Etau}
\tau:\  T_i&\mapsto T_i,\\
 X_j & \mapsto X_j
\end{split}\end{equation*}
for all $i = 0,\dots,n-1, j = 0,\dots,n$. Since the relations given in Section
1 are all invariant with respect to reversal of the order of generators, this does indeed define an antiautomorphism.

As any antiautomorphism, $\tau$ can be used to define a left action of $\H_n$ on
the $F$-dual $M^*=\hom_F(M,F)$ of a module $M \in \H_n \mod$ via $hf(m)=f(\tau(h)m)$. Denote this module by $M^\tau$.
As $\tau$ leaves parabolic subalgebras of $\H_n$ invariant,
it can also be used to define a duality on finite dimensional $\H_I$-modules
for any subset $I$ of $\{0,\dots,n-1\}$. If we think of representations in
terms of matrices, the $\tau$-dual corresponds to taking the transposes of the
representing matrices. Then we obtain another corollary of Lemma \ref{joe3}:

\begin{mycor}\label{c2}
For $M \in \H_I \mod$, there is a natural isomorphism 
$$(\ind_I^n M)^\tau \cong \ind_I^n ({^d}(M^\tau)).$$
\end{mycor}

\begin{proof} 
Recall first that $I \cap dI=I$, so the induction on the right actually makes sense. Then the proof is analogous to \cite{BK:01}, Lemma 2.13.
\end{proof}

\section{Formal Characters}

In this section  we investigate how far formal characters - mainly tuples of eigenvalues -, which uniquely determine irreducible modules for the affine Hecke algebra of type $A$, lead in type $B$. 
In the following, we will make heavy use of the following well-known lemma.

\begin{lemma}\label{fear}
For $F$-algebras $A$ and $B$, the irreducibles in $A \otimes_F B \mod$ are
exactly the outer tensor products $M \boxtimes N$ of irreducible $M \in A
\mod$, $ N \in B \mod$. Further, if $M \boxtimes N \cong M' \boxtimes N'$,
then $M \cong M'$ and $N \cong  N'$.
\end{lemma}

For $\ba=(a_0, a_1,\dots,a_n)\in F^{n+1}$, the one-dimensional $\P_n$-module
on which $X_i$ acts as the scalar $a_i$ for $0 \leq i \leq
n$ will also be denoted by $\ba$. If an eigenvalue occurs several
times this will be indicated by an exponent in parentheses. So $(a_0, a^{(n)})$
is the one-dimensional $\P_n$-module on which $X_0$ acts as $a_0$ and all $X_i$ for $i>0$ act as $a$.
Since $\P_n$ is commutative and $F$ is algebraically closed, the modules $(a_0,a_1, \dots, a_n)$ exhaust all
irreducibles in $\P_n \mod$. This follows from the fact that commuting
matrices can simultaneously be brought to upper triangular form.
 
For $M\in \P_n \mod$ and
 any $\ba \in F^{n+1}$,
let $M[{\ba}]$ be the largest submodule of $M$ all of whose 
composition factors are isomorphic to 
$\ba$, i.e.\ 
$$M[{\ba}]= \{m \in M \mid \exists k \in \ZZ_{\geq 0} \forall 0\leq i \leq n: (X_i
-a_i)^k m=0      \}$$
is the
simultaneous generalized $(a_0, a_1,\dots,a_n)$-eigenspace for
 $X_0,X_1,\dots, X_n$.

\begin{lemma}\label{decomp1}
For any $M \in \P_n \mod$ we have
$M = \oplus_{\ba\in F^n} M[\ba]$ as a $\P_n$-module.
\end{lemma}

\begin{proof}
As an  $F[X_i]$-module we can decompose $M$ into the direct sum of
generalized eigenspaces since $F$ is algebraically closed and we have a Jordan
normal form. Since $\P_n$ is commutative this decomposition is respected by
all the other $X_j$ and the assertion follows by induction. 
\end{proof}

We define the {\em formal character} of $M \in \H_n \mod$ by:
\begin{equation}
\ch M := [\res^{\H_n}_{\P_n} M] \in K(\P_n\mod).
\end{equation}
Exactness of the functor $\res^{\H_n}_{\P_n}$ implies that 
$\ch$ induces a homomorphism
$$
\ch:K(\H_n\mod) \rightarrow K(\P_n \mod)
$$
between the corresponding Grothendieck groups. 
For $M \in \H_I \mod$, the definition is modified in the obvious way.
Note that if we expand 
$$\ch M = \sum_{\ba \in F^{n+1}} r_{\ba} 
[(a_0,a_1,\dots,a_n)]$$ in terms of the basis for $K(\P_n\mod)$ given by the irreducibles,
 the coefficient $r_{\ba}$ is exactly the dimension of the generalized simultaneous 
$\ba$-eigenspace $M[\ba]$ of $X_0, \dots, X_n$.

We can explicitly compute formal characters of induced modules.

\begin{lemma}
\label{LChari}
Let $\ba = (a_0,a_1,\dots,a_n)\in F^{n+1}$. Then
$$
\ch \ind_{\P_n}^{\H_n}\ba
=\underset{\substack{u \in S_n \\ \underline{\eps} 
\in \{1, -1 \}^n }}{\sum}
[(b_0(u,\underline{\eps}), a_{u^{-1}(1)}^{\eps_1}, \dots, a_{u^{-1}(n)}^{\eps_n})]
$$ where $b_0(u,\underline{\eps}) := a_0 \underset{\substack{j \in \{1, \dots,n\} \\
  \eps_j=-1}}{\prod} a_{u^{-1}(j)}$.
\end{lemma}
\begin{proof}
This follows directly from the Mackey Theorem with $I =  J = \emptyset$.
\end{proof}

\begin{lemma}{\rm (``Shuffle Lemma'')} \label{shuffle}

Let $n=m+k$, and let $M\in \H_m\mod$, $K\in\H_k^A\mod$. Assume 
$$
\ch M=\sum_{{\ba}\in F^{m+1}}r_{\ba}[(a_0,\dots,a_m)],\qquad
\ch K=\sum_{{\bb}\in F^k}t_{\bb}[(b_1,\dots, b_k))].
$$
Then 
$$
\ch \ind_{m,k}^n M\boxtimes K = \sum_{{\ba}\in F^m}\sum_{{\bb}\in F^k}
r_{\ba} t_{\bb}(\sum_{\bc} (c_0,c_1, \dots,c_n)),
$$
where the last sum is over all $\bc=(c_0,\dots,c_n) \in F^{n+1}$ such that  $(c_1,\dots,c_n)$ is obtained obtained
by shuffling $(a_1, \dots, a_n)$ and $\bb^{\underline{\eps}}:=(b_1^{\eps_1}, \dots, b_k^{\eps_k})$ for $\underline{\eps} =(\eps_1, \dots , \eps_k) 
\in \{1, -1 \}^k$,
which means there exist $1\leq u_1<\dots<u_m\leq n$ 
such that $(c_{u_1},\dots,c_{u_m})=(a_1,\dots,a_m)$,
$(c_1,\dots,\widehat {c}_{u_1},\dots,\widehat{c}_{u_m},\dots,c_n)=
\bb^{\underline{\eps}}$ and $c_0 = a_0 \underset{\substack{j \in \{1, \dots,k\} \\
  \eps_j=-1}}{\prod} b_{j}$.
\end{lemma}

\begin{proof}
This also follows directly from the Mackey Theorem setting $I = \emptyset$ and
$ J = (m,k)$.
\end{proof}

At this point, it is convenient to introduce the left coset representatives of
the two maximal parabolic subgroups that will be most important in the
following. 

The set $D_{(n-1,1)}$  of distinguished left coset representatives of $W_{n-1}^\fn$ in $W_n^\fn$ consists of $1$, $s_{j,n-1}$ for $0 \leq j \leq n-1$ and
$s_{j,0,n-1}$ for $1 \leq j \leq n-1$, see e.g.\ \cite{DJ:92}.

The set $D_{(0,n)}$ can easily be checked to consist of all elements of the form$$s_{\underbar j}:= s_{j_r,0}s_{j_{r-1},0} \cdots s_{j_1,0}$$
for subsets $\underbar j=(j_1 > j_2> \dots > j_r) \subseteq \{0, \dots,n-1\}$.

To illustrate the way of computing formal characters using the Shuffle Lemma,
we will now give two examples of induced modules from parabolic subalgebras
of the types described above.
\begin{example}
 First we take the irreducible module for $\H_{1,1}$ which is isomorphic to the
outer tensor product of the one-dimensional $\H_1$-module $L(a_0,p^2)$ on which $T_0$ acts as $p$, $X_0$ as $a_0$ and $X_1$ as $p^2$ and the one-dimensional $\H_1^A \cong F[X_2^{\pm1}]$-module $(c)$ on which $X_2$ acts by multiplication with $c$.
The set $D_{\emptyset,(1,1)}$, which  is simply $D_{(1,1)}$ can be ordered as $(1,s_1,s_0s_1,s_1s_0s_1)$ and the Shuffle Lemma yields

\begin{equation*}\begin{split}
\ch \ind_{\H_{1,1}}^{\H_2} L(a_0,p^2) \boxtimes (c)
&=[(a_0,p^2,c)] + [(a_0,c,p^2)]\\
& \quad + [(a_0c,c^{-1},p^2)] + [(a_0c,p^2,c^{-1})] .
\end{split}\end{equation*}

As a second example, we consider $\ind_{\P_0 \otimes \H_3^A}^{\H_3}
(a_0) \boxtimes L^A(-q^{-2},-1,-q^2)$, where $L^A(-q^{-2},-1,-q^2)$ denotes the
one-dimensional $\H_3^A$-module, on which $T_1$ and $T_2$ act by
multiplication with $q$, and  $X_1,X_2,X_3$ act as $-q^{-2},-1,-q^2$
respectively.
Ordering the set $D_{\emptyset,(0,3)}=D_{(0,3)}$ as
$$(1,s_0,s_1s_0,s_2s_1s_0,s_0s_1s_0,s_0s_2s_1s_0,s_1s_0s_2s_1s_0, s_0s_1s_0s_2s_1s_0),$$ we obtain
\begin{equation*}\begin{split}
\ch \ind_{\P_0 \otimes \H_3^A}^{\H_3} &(a_0) \boxtimes L^A(-q^{-2},-1,-q^2)\\
&=[(a_0,-q^{-2},-1,-q^2)] + [(-a_0q^{-2},-q^2,-1,-q^2)]\\
& \quad + [(-a_0q^{-2},-1, -q^2,-q^2)] + [(-a_0q^{-2},-1, -q^2,-q^2)] \\
&  \quad + [(a_0q^{-2},-1, -q^2,-q^2)]+ [(a_0q^{-2},-1, -q^2,-q^2)]\\
& \quad + [(a_0q^{-2}, -q^2,-1,-q^2)]+ [(-a_0, -q^{-2},-1,-q^2)].\\
\end{split}\end{equation*}
\end{example}

Recall that the center $Z(\H)$ of $\H$ consists of
those Laurent polynomials $f$ which are invariant under the action of $W_n^\fn$.
Given $\ba \in F^{n+1}$, we associate the
{\em central character}
$$
\chi_{\ba}:Z(\H_n) \rightarrow F,
\quad
f(X_0,\dots,X_n) \mapsto f(a_0,\dots,a_n),
$$
so central characters are simply all algebra homomorphisms from $Z(\H_n)$ to $F$.
Consider the left action of $W_n^\fn$ on $F^{n+1}$ induced by the action on $\P_n$, i.e.\ $s_i(a_0,\dots,a_n) = (a_0,\dots,a_{i+1},a_i,\dots,a_n)$ for $1 \leq i \leq n-1$ and 
$s_0(a_0,\dots,a_n) = (a_0a_1,a_1^{-1},a_2,\dots,a_n)$.
Writing $\ba \sim \bb$ if $\ba$ and $\bb$ lie in the same orbit with respect to this action, we get:

\begin{lemma}\label{ccc} For $\ba,\bb \in F^{n+1}$, 
$\chi_{\ba} = \chi_{\bb}$ if and only
if $\ba \sim \bb$.
\end{lemma}

Thus the central characters of $\H_n$ are actually labeled by the set 
$F^{n+1}/{\sim}$ of $W_n^\fn$-orbits on $F_{n+1}$ and we set, for $\ga \in F^{n+1}/{\sim}$
$$
\chi_\ga:=\chi_\ba
$$ 
for any $\ba\in\ga$.

Accordingly, for $M \in \H_n \mod$ we define
$$
M[\gamma] = \{v \in M\:|\: (z - \chi_\gamma(z))^k v = 0
\hbox{ for all $z \in Z(\H_n)$ and $k \gg 0$}\}.
$$
It is important that this is an $\H_n$-submodule of $M$. Indeed, $(z - \chi_\gamma(z))^k v = 0$ implies 
$(z - \chi_\gamma(z))^k hv = h(z - \chi_\gamma(z))^k v = 0$ for all $h\in \H_n$.
Now, for $\ba \in F^{n+1}$ with $\ba \in \gamma$,
$Z(\H_n)$ acts on the $\P_n$-submodule 
$M[\ba]$ via the central
character $\chi_\gamma$. Hence
$$
M[\gamma] = \bigoplus_{\ba\in\gamma} M[\ba],
$$
as a $\P_n$-module.
So, recalling Lemma \ref{decomp1}, we see:
\begin{lemma}
\label{L070900}
Any $M\in\H_n\mod$ decomposes as
$$
M = \bigoplus_{\gamma \in F^{n+1}/{\sim}} M[\gamma]
$$
as an $\H_n$-module.
\end{lemma}

Thus the $\chi_{\gamma}$, $\gamma \in F^{n+1}/{\sim}$,
exhaust the possible central characters that can arise
in a finite dimensional $\H_n$-module, while
Lemma \ref{LChari} shows that every such central character
does arise in some finite dimensional $\H_n$-module.

For $\gamma\in F^{n+1}/\sim$, we define $\H_n[\gamma]\mod$ to be the full subcategory of
$\H_n\mod$ consisting of all modules $M$ with $M[\gamma] = M$.
Lemma \ref{L070900} in fact yields an equivalence of categories
\begin{equation}\label{block1}
{\H_n} \mod \cong
\bigoplus_{\gamma \in F^{n+1}/{\sim}} {\H_n}[\ga]\mod.
\end{equation}

In the following, we will call $\H_n[\gamma]\mod$ the block of $\H_n \mod$ corresponding to $\ga$.

Observe that none of this actually uses which affine Hecke algebra we work with, 
thus the analogous statements hold for any parabolic subalgebra of $\H_n$ 
as well as of $\H_n^R$, in particular we have the same notions of formal characters, central characters and blocks for $\H_n^A$.


We now need a type $A$ result for a special module that will be important in the following section.
It is due to Kato \cite{K:81}, but for convenience, we include a proof. 
Denote by $ L^A(a^{(n)}):= \ind_{\R_n}^{ \H^A_n}( a^{(n)})$ the principal
series module in type $A$.

\begin{lemma}\label{TKato} Let $a\in F$
and $I = (0, m_1, \dots, m_r) \subseteq \{0, \dots, n-1\}$. 
\begin{enumerate}
\item[(i)]
$L^A(a^{(n)})$ is irreducible, and it is the only irreducible module in its
block.
\item[(ii)]
All composition factors of 
$\res^{ \H^A_n}_{ \H^A_{I}} L^A(a^{(n)})$ 
are isomorphic to $$L^A(a^{(m_1)})\boxtimes \dots\boxtimes L^A(a^{(m_r)}),$$
and $\soc \res^{ \H^A_n}_{ \H^A_{I}} L^A(a^{(n)})$
is irreducible.
\item[(iii)]
$\soc \res^{ \H^A_n}_{ \H^A_{n-1}} L^A(a^{(n)}) \cong L^A(a^{(n-1)})$.
\item[(iv)] The size of any Jordan block of $X_n$ on $L^A(a^{(n)})$ is $n$.
\end{enumerate}
\end{lemma}

\begin{proof}
We first show that the $a^{(n-1)}$-eigenspace for $X_1, \dots, X_{n-1}$ is contained in $1 \otimes (a^{(n)})$, which is, of course, also contained in  the $a$-eigenspace of $X_n$.
This is certainly true for $n=1$, since in this case $L^A(a) =(a)$ as $\H_1^A = \R_1$.
So assume that the statement holds for $n-1$, i.e.\ the $a^{(n-2)}$-eigenspace for $X_1, \dots, X_{n-2}$ in $L^A(a^{(n-1)})$ is contained in $1 \otimes (a^{(n-1)})$. Let $u \in L^A(a^{(n)})$ be an $a^{(n-1)}$-eigenvector for $X_1, \dots, X_{n-1}$. As $L^A(a^{(n)}) \cong \ind_{\H_{n-1,1}^A}^{ \H^A_n}L^A( a^{(n-1)})\boxtimes (a)$, $u$  can be written as
\begin{equation*}\begin{split}
u= \sum_{j=1}^{n-1}T_{j,n-1} \otimes u_j + 1 \otimes u_0
\end{split}\end{equation*}
with $u_j \in L^A(a^{(n-1)}) \boxtimes (a)$ for $0 \leq j \leq n-1$. 
Now suppose $l$ is minimal such that $u_l \neq 0$, then, if $l \le n-1$
\begin{equation*}\begin{split}
(X&_{l+1}-a)(T_{l,n-1} \otimes u_l + T_{l+1,n-1} \otimes u_{l+1})\\
& = T_{l,n-1}  (X_{l}- a) \otimes u_{l}  + (q-q^{-1})T_{l+1,n-1} X_{n}
\otimes u_{l}\\ &\quad+ T_{l+1,n-1} (X_n-a) \otimes u_{l+1}\\ &\quad+ \text{terms in $T_{j,n-1} \otimes L^A( a^{(n-1)})\boxtimes (a)$ for $j \geq l+2$ and $1 \otimes L^A( a^{(n-1)})$}.
\end{split}\end{equation*}
Since $X_n$ acts as $a$ on $L^A( a^{(n-1)})\boxtimes (a)$, we see from the coefficient of $T_{l+1,n-1}$ that $u_l$ has to be zero, a contradiction. 
If $l =n-1$, then for $k \leq n-2$
\begin{equation*}\begin{split}
(X_{k}&-a)(T_{n-1} \otimes u_{n-1} + 1 \otimes u_{0})\\
&= T_{n-1}  (X_{k}- a) \otimes u_{n-1}  + 1 \otimes (X_{k}-a) u_{0},
\end{split}\end{equation*}
whence, $u_{n-1}, u_0$ have to be contained in $1 \otimes (a^{(n-1)})\boxtimes (a) \subseteq L^A( a^{(n-1)})\boxtimes (a)$ by the inductive hypothesis. But 
\begin{equation*}\begin{split}
(X_{n-1}&-a)(T_{n-1} \otimes u_{n-1} + 1 \otimes u_{0})\\
&=  T_{n-1}  (X_{n}- a) \otimes u_{n-1} \\& \quad- (q-q^{-1}) 1 \otimes X_n u_{n-1} + 1 \otimes (X_{n-1}-a) u_{0},
\end{split}\end{equation*}
requiring $(X_{n-1}-a) u_{0} =(q-q^{-1}) 1 \otimes X_n u_{n-1} $, but\mbox{$(X_{n-1}-a) u_{0}=0$} as we have just seen. 
Hence, $l$ does not exist, so the $a^{(n)}$-eigenspace is contained in $1 \otimes L^A( a^{(n-1)})\boxtimes (a)$ and therefore, by induction, in $1 \otimes (a^{(n)})$.

Now any nontrivial submodule of $L^A(a^{(n)})$ has to contain a simultaneous eigenvector for $X_1, \dots , X_n$ which can only be an $a^{(n)}$-eigenvector as the formal character of $L^A( a^{(n)})$ is $n![(a^{(n)})]$. Therefore any nontrivial submodule contains $1 \otimes (a^{(n)})$, generating the whole of
$L^A( a^{(n)})$. This proves (i).

The assertion that all composition factors of 
$\res^{ \H^A_n}_{ \H^A_{I}} L^A(a^{(n)})$ 
are isomorphic to \mbox{$L^A(a^{(m_1)})\boxtimes \dots\boxtimes L^A(a^{(m_r)})$} follows immediately from the Mackey Theorem and the irreducibility of the latter module by (i) and Lemma \ref{fear}. Therefore, the socle of $\res^{ \H^A_n}_{ \H^A_{I}} L^A(a^{(n)})$ consists of a number of copies of $L^A(a^{(m_1)})\boxtimes \dots\boxtimes L^A(a^{(m_r)})$. But any constituent of the socle contains a simultaneous eigenvector for $\R_n$, but there is only one such, up to a scalar multiple, hence the socle is irreducible, completing the proof of (ii).

By (ii) the socle of $\res^{ \H^A_n}_{ \H^A_{n-1,1}} L^A(a^{(n)})$ is isomorphic to $L^A(a^{(n-1)}) \boxtimes (a)$, hence the socle of  $\res^{ \H^A_n}_{ \H^A_{n-1}} L^A(a^{(n)})$ certainly contains a copy of $L^A(a^{(n-1)})$. But again, every constituent of $\soc \res^{ \H^A_n}_{ \H^A_{n-1}} L^A(a^{(n)})$ contains a simultaneous eigenvector of $\R_{n-1}$, which also is contained in $1 \otimes (a^{(n)})$, hence 
$\soc \res^{ \H^A_n}_{ \H^A_{n-1}} L^A(a^{(n)}) \cong L^A(a^{(n-1)})$.

For the proof of (iv), assume inductively that the size of any Jordan block of $X_{n-1}$ on $L^A(a^{(n-1)})$ is $n-1$ which we can do since the statement certainly holds for the $\H_1^A$-module $(a)$. Hence there are $(n-2)!$ elements in $L^A(a^{(n-1)})$, each generating a Jordan block of size $n-1$ for $X_{n-1}$ and similarly for $L^A(a^{(n-1)})\boxtimes (a)$.
Now for such an element $v$ and any $1 \leq l \leq n-1$
\begin{equation*}
(X_{n}-a)T_{l,n-1} \otimes v
=T_{l,n-1}  (X_{n-1}- a) \otimes v  + (q-q^{-1})T_{l,n-2} (X_{n})
\otimes v
\end{equation*}
whence 
\begin{equation*}
(X_{n}-a)^{n-1}T_{l,n-1} \otimes v
= (q-q^{-1})T_{l,n-2} \otimes (X_{n})  (X_{n-1}- a)^{n-2} v
\end{equation*}
which is nonzero but annihilated by $(X_{n}-a)$. So we have \mbox{$(n-1)(n-2)! = (n-1)!$} Jordan blocks of size $n$ for $X_n$, exhausting $L^A(a^{(n)})$.
\end{proof}

In type $B$, we need certain restrictions on the eigenvalue $a$ to prove a similar statement.

\begin{lemma}\label{Kato}
Let $a\in F \setminus \{\pm1\}$. Then
\begin{enumerate}
\item $\ind_{\P_n}^{\H_n} (a_0, a^{(n)})$ has an irreducible cosocle and 
\item this cosocle, denoted by $L(a_0,a^{(n)})$, is the only irreducible module in $\H_n \mod$ containing $[(a_0, a^{(n)})]$ as a formal character.
\end{enumerate}

\end{lemma}

\begin{proof}
Note that
\begin{equation*}\begin{split}
\ind_{\P_n}^{\H_n} (a_0, a^{(n)})& \cong \ind_{\P_0 \otimes \H^A_n}^{\H_n}
\ind_{\P_n}^{\P_0 \otimes \H^A_n}(a_0, a^{(n)}) \\
& \cong \ind_{\P_0 \otimes \H^A_n}^{\H_n}(a_0)\boxtimes L^A(a^{(n)}).
\end{split}\end{equation*}

Now, applying the Shuffle Lemma, we see that in fact, the only formal characters of
the form $[(a_0, a^{(n)})]$ in $\ind_{\P_0 \otimes \H^A_n}^{\H_n}(a_0)\boxtimes
L^A(a^{(n)})$ stem from the coset representative $1$, from which it follows that such formal character values can occur only in the cosocle of the induced module. Knowing this, (ii) will follow as soon as we have established (i). But from Frobenius reciprocity (\ref{frob}), we have
\begin{equation*}\begin{split}
\hom_{\H_n}&(\ind_{\P_n}^{\H_n} (a_0, a^{(n)}),\cosoc \ind_{\P_n}^{\H_n}
(a_0, a^{(n)}))\\ &\cong \hom_{\P_0 \otimes \H_n^A}
((a_0) \boxtimes L^A(a^{(n)}),\res^{\H_n}_{\P_0 \otimes \H_n^A} \cosoc \ind_{\P_n}^{\H_n}
(a_0, a^{(n)}) )\\& = F,
\end{split}\end{equation*}
proving (i). 
\end{proof}

\section{``Kashiwara operators''}

In this section we will, in analogy to the type $A$ situation, define maps between the sets of isomorphism classes of irreducibles in $\H_n \mod$ and $\H_{n-1} \mod$. 

Let $M \in \H_n\mod$ and $a \in F$.
Define $\De_a M$ to be the generalized $a$-eigenspace of $X_n$ in $M$, i.e.\
$$
\De_a M = \bigoplus_{\ba \in F^{n+1},\ a_n = a} M[\ba].
$$

As $X_n$ is central in the parabolic subalgebra
$\H_{n-1,1} \cong \H_{n-1} \otimes F[X_n^{\pm1}]$ of $\H_n$, 
$\De_a M$ is an $\H_{n-1,1}$-submodule of $M$.
Thus $\De_a$ defines a functor
\begin{equation*}
\De_a:\H_n\mod \rightarrow \H_{n-1,1}\mod,
\end{equation*}
which, on morphisms, is simply restriction. This functor is exact as the composite of restriction to a subalgebra and then taking a direct summand.
Analogously, for $m \geq 0$,  $\De_{a^{(m)}}$ denotes the functor $\H_{n}\mod \rightarrow \H_{n-m,m}\mod$ that maps $M$ to simultaneous generalized $a$-eigenspace of $X_{n-m+1},\dots,X_n$.

By Lemmas \ref{fear} and \ref{TKato}, $\De_{a^{(m)}} M$ is the largest submodule of $\res^n_{n-m,m} M$ 
all of whose composition factors are of the form $N\boxtimes L^A(a^{(m)})$
for some irreducible \mbox{$N \in \H_{n-m}\mod$} and is indeed a direct summand of $\res^n_{n-m,m} M$. 

 \begin{lemma}\label{E090900_1}
For $N\in \H_{n-m}\mod$, $M\in\H_n\mod$, there is a functorial isomorphism
\begin{equation*}
\hom_{\H_{n-m,m}}(N \boxtimes L^A(a^{(m)}),\De_{a^{(m)}} M)
\cong \hom_{\H_n}(\ind_{n-m,m}^n N \boxtimes L^A(a^{(m)}), M).
\end{equation*}
\end{lemma} 

\begin{proof}
By Lemma \ref{fear}, all composition factors of $\res_{n-m,m}^n M$ are of the
form $K \boxtimes L$ for irreducible $K \in \H_{n-m}\mod$ and $L \in \H_m^A
\mod$. An injection of the irreducible $N \boxtimes L^A(a^{(m)})$ into
$\res_{n-m,m}^n M$ can only map onto a submodule that is isomorphic to $N
\boxtimes L^A(a^{(m)})$ and all composition factors with this property are
contained in $\De_{a^{(m)}} M$. Since $\De_{a^{(m)}} M$ is a direct summand of
$\res_{n-m,m}^n M$, the assertion follows from Frobenius reciprocity (\ref{frob}).
\end{proof}

The following is immediate from the definition:
\begin{lemma}\label{L310800_2} 
Let $M\in\H_n\mod$ with 
$$\ch M=\sum_{{\ba}\in F^{n+1}}r_{\ba}[(a_0,\dots,a_n)].$$
Then we have 
$$\ch \De_{a^{(m)}} M=
\sum_{\bb} r_{\bb}[(b_0,\dots,b_n)],
$$
summing over all $\bb \in F^{n+1}$ with
$b_{n-m+1} =\dots=b_n = a$.
\end{lemma}

Now for $a \in F$ and $M \in \H_n\mod$, 
define
\begin{equation*}
\label{E280900} 
\eps_a(M)=\max\{m\geq 0\:|\: \Delta_{a^{(m)}} M\neq 0\}. 
\end{equation*} 
By Lemma \ref{L310800_2}, $\eps_a(M)$ is simply the length of the `longest
$a$-tail' in $\ch M$.

\begin{lemma}\label{sick} 
Let $M\in \H_n\mod$ be irreducible, $a\in F$, $\eps=\eps_a(M)$.
If \mbox{$N \boxtimes L^A(a^{(m)})$} is an irreducible
submodule of $\De_{a^{(m)}} M$ for some $0 < m \leq \eps$ and some irreducible $N \in \H_{n-m}\mod$,
then $\eps_a(N) =\eps-m$.  
\end{lemma} 

\begin{proof} 
Analogous to \cite{BK:01}, Lemma 5.2.
\end{proof}

\begin{lemma}\label{L010900}
Let $m \geq 0$, $a \in F \setminus \{\pm 1\}$ and $N\in \H_{n-m} \mod$ be irreducible 
with $\eps_a(N) = 0$.
Set $M = \ind_{n-m,m}^{n} N \boxtimes L^A(a^{(m)})$.
Then 
\begin{enumerate}
\item[(i)] $\De_{a^{(m)}} M \cong N \boxtimes L^A(a^{(m)})$ ;
\item[(ii)]$\cosoc M$ is irreducible with $\eps_a(\cosoc M) = m$;
\item[(iii)] All other composition factors $L$ of $M$ have $\eps_a(L) < m$.
\end{enumerate}
\end{lemma}

\begin{proof} 
Analogous to \cite{BK:01}, Lemma 5.3.
\end{proof}

\begin{lemma}
\label{L310800_1}
Let $M\in \H_n \mod$ be irreducible, $a \in F \setminus \{\pm 1\}$, $\eps=\eps_a(M)$ and $0 \leq m \leq \eps$. Then
\begin{enumerate}
\item[(i)] 
$\De_{a^{(\eps)}} M \cong N \boxtimes L^A(a^{(\eps)})$ for some
irreducible  $N \in \H_{n-\eps}\mod$ with \mbox{$\eps_a(N) = 0,$}
\item[(ii)]  $\soc \De_{a^{(m)}} M \cong L \boxtimes L^A(a^{(m)})$ for some irreducible $L \in \H_{n-m} \mod$ with 
$\eps_a(L) = \eps_a(M) - m$.
\end{enumerate}
In particular, if $\pm 1$ do not occur as eigenvalues of $X_n$ on $M$, $\soc \res^n_{n-1,1} M$ is multiplicity-free.
\end{lemma}

\begin{proof}
(i) Pick an irreducible submodule $N\boxtimes L^A(a^{(\eps)})$ of
$\De_{a^{(\eps)}} M$, cf.\ remark before Lemma \ref{E090900_1}.
Then $\eps_a(N) = 0$ by Lemma~\ref{sick}. 
By Lemma \ref{E090900_1}
$$\hom_{\H_{n-\eps,\eps}}(N \boxtimes L^A(a^{(\eps)}),\De_{a^{(\eps)}} M)
\cong \hom_{\H_n}(\ind_{n-\eps,\eps}^n N \boxtimes L^A(a^{(\eps)}), M) \neq 0,$$
thus $M$, being irreducible, is a quotient of  $\ind_{n-\eps,\eps}^n N\boxtimes L^A(a^{(\eps)})$.
Hence, $\De_{a^{(\eps)}} M$ is a quotient of
$\De_{a^{(\eps)}} \ind_{n-\eps,\eps}^n N\boxtimes L^A(a^{(\eps)})$.
But, by Lemma~\ref{L010900}(i), this is irreducible and isomorphic to 
$N \boxtimes L^A(a^{(\eps)})$, proving (i).

(ii) For every constituent $L\boxtimes L^A(a^{(m)})$ of $\soc \De_{a^{(m)}} M$, Lemma \ref{sick} tells us that $\eps_a(L)=
\eps-m$, hence $\De_{a^{(\eps -m)}} L \boxtimes L^A(a^{(m)})$ is a non-trivial submodule 
of $\res^{n-\eps,\eps}_{n-\eps,\eps-m,m}\De_{a^{(\eps)}} M$. 
By (i), $\De_{a^{(\eps)}} M$  is irreducible of the form $N\boxtimes L^A(a^{(\eps)})$, so Lemma \ref{TKato}(ii) implies that 
$$ \soc \res^{n-\eps,\eps}_{n-\eps,\eps-m,m}\De_{a^\eps} M \cong N\boxtimes L^A(a^{(\eps-m)})\boxtimes L^A(a^{(m)}).$$ 
Therefore there can only be one such constituent in $\soc\De_{a^{(m)}} M$ and we're done.
\end{proof}

Defining  
\begin{equation}
\label{EE}
e_a:=\res^{n-1,1}_{n-1}\circ \Delta_a : \H_n\mod \rightarrow \H_{n-1}\mod
\end{equation}
and analogously

$$f_a= \ind_{n-1,1}^n - \boxtimes (a) : \H_{n-1}\mod \rightarrow \H_n\mod $$
we obtain the following corollary.

\begin{mycor}\label{karg}
For $a \in F \setminus \{\pm 1\}$ and an irreducible $M \in \H_n\mod$ with \mbox{$\eps_a(M) > 0$,} 
$\soc e_a M$ is irreducible, and $\eps_a(\soc e_a M)=\eps_a(M)-1$. 
\end{mycor}

\begin{proof}
Choose a constituent $L$ of $\soc e_a M$. 
The central element $Z:=X_0^2X_1\dots X_n$ of $\H_n$ acts as a scalar on the 
whole of $M$ by Schur's Lemma, and similarly 
the central element \mbox{$Z':=X_0^2X_1\dots X_{n-1}$} of $\H_{n-1}$ 
acts as a scalar on $L$. Hence $X_n=Z'^{-1}Z$ acts on 
$L$ as a scalar, too. The scalar must be $a$, so  $L$ contributes a 
constituent \mbox{$L\boxtimes (a)$} to 
$ \soc \De_a M$. By Lemma \ref{L310800_1} (ii) this is irreducible and satisfies \mbox{$\eps_a(\soc \De_a M)=\eps_a(M)-1$,}  so $\soc e_a M$ must also be irreducible and isomorphic to $L$ and $\eps_a(\soc e_a M)=\eps_a(M)-1$.
\end{proof} 

The following lemma provides a recipe for an inductive construction of irreducible modules in $\H_n \mod$ from irreducibles in $\H_{n-1}\mod$.

\begin{lemma}
\label{L090900_1}
Let $m\geq 0$, 
$a \in F \setminus \{\pm 1\}$, let $N\in\H_n\mod$ be irreducible and set 
$M=\ind_{n,m}^{n+m}(N \boxtimes L^A(a^{(m)})).$
Then, $\cosoc M$ is irreducible with $\eps_a(\cosoc M)=
\eps_a(N)+m$, and all other composition factors $L$ of 
$M$ have $\eps_a(L)<\eps_a(N)+m$. 
\end{lemma}

\begin{proof}
Analogous to \cite{BK:01}, Lemma 5.5.
\end{proof}

We can now define the desired Kashiwara type operators.
Let $M$ be an irreducible module in ${\H_n}\mod$. Define
\begin{equation}\label{tildeef}
\tilde e_a M:=\soc e_a M,\quad
\tilde f_a M:=\cosoc f_a M = \cosoc \ind_{n,1}^{n+1} M \boxtimes (a),
\end{equation}
For $a \neq \pm 1$, $\tilde f_a M$ is irreducible by 
Lemma \ref{L090900_1} and $\tilde e_a M$ is irreducible or $0$ by 
Corollary \ref{karg}, hence the functors induce maps between the set of isomorphism classes of irreducibles in ${\H_n}\mod$ and ${\H_{n-1}}\mod$.
Observe that Corollary \ref{karg} implies 
\begin{equation}
\label{E300900}
\eps_a(M)=\max\{m\geq 0\:|\: \tilde e_a^m M\neq 0\}
\end{equation}
and, by Lemma \ref{L090900_1}, we have
\begin{equation}
\label{hi}
\eps_a(\tilde f_a M)=\eps_a (M) + 1.
\end{equation}

\begin{lemma}\label{sunny}
Let $M\in{\H_n}\mod$ be irreducible, $a \in F \setminus\{\pm1\}$ and $m \geq 0$.
\begin{enumerate}
\item[(i)]
$\soc \De_{a^{(m)}} M \cong (\tilde e_a^m M)\boxtimes L^A(a^{(m)}).$
\item[(ii)]
$\cosoc \ind_{n,m}^{n+m} M \boxtimes L^A(a^{(m)}) \cong \tilde f_a^m M.$
\end{enumerate}
\end{lemma} 
\begin{proof} Analogous to \cite{BK:01}, Lemma 5.9.\end{proof}

\begin{lemma}
\label{L290900_2}
Let $M\in {\H_n}\mod$ and $N \in {\H_{n-1}}\mod$ be irreducible modules
and $a  \in F \setminus \{\pm1\}$. Then, 
$\tilde e_a M\cong N$ if and only if $\tilde f_a N\cong M$. 
\end{lemma}

\begin{proof}
Analogous to \cite{BK:01}, Lemma 5.10.
\end{proof}

We immediately get the following corollary:

\begin{mycor}\label{imm}
Let $M, N\in{\H_n}\mod$ be irreducible and \mbox{$a \in F \setminus \{\pm1\}$.}
Then
\begin{enumerate}
\item $\tilde e_a \tilde f_a M \cong M$ and, if $\eps_a(M) > 0$, $\tilde f_a \tilde e_a M \cong M$;
\item $\tilde f_a M \cong \tilde f_a N$ if and only if
$M \cong N$ and, if $\eps_a(M),\eps_a(N) > 0$,
$\tilde e_a M \cong \tilde e_a N$ if and only if $M \cong N$.
\end{enumerate}
\end{mycor}

We define the graph $\Ga$ as the graph whose vertices correspond to
isomorphism classes of irreducible modules in $\underset{n\geq 0}{\bigoplus}\H_n \mod$, where there is an edge \mbox{$[N] \overset{a}{\rightarrow} [M]$} if and only if $M \cong \tilde f_a N$. 
If an irreducible  $M \in \H_n \mod$ is isomorphic to $ \tilde f_{a_n} \cdots \tilde f_{a_1} (a_0)$, we write $M \cong L(a_0,a_1, \dots,a_n)$. 
Note that these definitions make sense even in the case where $a \in
\{\pm1\}$. Even though in this case $\tilde f_a$ does not necessarily produce
irreducible modules, as we will see in Section 5, we just don't draw any edges
from $N$ if $\cosoc \ind_{n-1,1}^n N \boxtimes (a)$ is not irreducible. This has the drawback that not every vertex in $\Ga$ is connected to a module for $\H_0$, but it enables us to use the (by this definition unique) labeling $L(\ba)$ in the cases where $1$ or $-1$ occur in $\ba$, but the corresponding cosocle of the induced modules are irreducible. 

We will denote the full subcategory of $\H_n \mod$ consisting of those modules on which $(X_i \pm 1)$ acts invertibly for all $1 \leq i \leq n$ by $\Rep_{\neq \pm 1} \H_n$. Then we obtain the following:

\begin{theorem}\label{inj}
The map
$\ch:K(\Rep_{\neq \pm 1}{\H_n}) \rightarrow K(\P_n \mod)$
is injective.
\end{theorem}
\begin{proof}
Analogous to \cite{BK:01}, Theorem 5.12.
\end{proof}

\begin{mycor}\label{selfdual}
If $L$ is an irreducible module in $\Rep_{\neq \pm 1}{\H_n}$,
then $L \cong L^\tau$.
\end{mycor}

\begin{proof}
Since $\tau(X_i) = X_i$, $\tau$ leaves characters invariant.
Hence it leaves irreducibles invariant since they are determined
up to isomorphism by their character according to the theorem.
\end{proof}

We now give three interpretations of the functions 
$\eps_a$.

\begin{theorem}
\label{T280900}
Let $a \in F \setminus \{\pm 1\}$ and $M$ be an irreducible module in ${\H_n}\mod$.
Then
\begin{enumerate}
\item[(i)] $[e_a M] = \eps_a(M) [\tilde e_a M] + \sum c_r [N_r]$
where the $N_r$ are irreducible modules with 
$\eps_a(N_r) < \eps_a(\tilde e_a M)$;
\item[(ii)] $\eps_a(M)$ is the maximal size of a Jordan block of
$X_n$ on $M$ with eigenvalue $a$;
\item[(iii)] The algebra $\End_{\H_{n-1}}(e_a M)$ is isomorphic to the algebra of truncated polynomials $F[x]/(x^{\eps_a(M)})$.
\end{enumerate}
\end{theorem}

\begin{proof}
Let $\eps = \eps_a(M)$ and $N = \tilde e_a^\eps M$.

(i) 
By Lemma~\ref{L310800_1} and Frobenius reciprocity, there is
a short exact sequence
$$
0 \longrightarrow K \longrightarrow \ind_{n-\eps,\eps}^n N \boxtimes 
L^A(a^{(\eps)}) \longrightarrow M \longrightarrow 0,
$$
where all composition factors $L$ of $K$ have
$\eps_a(L) < \eps$ by Lemma~\ref{L010900}(iii).
Applying the exact functor $\De_a$, we obtain the exact sequence
$$
0 \longrightarrow \De_a K \longrightarrow 
\De_a \ind_{n-\eps,\eps} N \boxtimes L^A(a^{(\eps)}) \longrightarrow \De_a M
\longrightarrow 0.
$$
By the Mackey Theorem
\begin{equation*}\begin{split}
 [ \res_{n-1,1}^n \ind_{n-\eps,\eps}^{n} N &\boxtimes L^A(a^{(\eps)})] = [\ind_{n-\eps,\eps-1,1}^{n-1,1} N\boxtimes \res_{\H^A_{\eps-1,1}}^{\H^A_\eps} L^A(a^{(\eps)})] \\
  &+ [\ind_{n-\eps-1,\eps,1}^{n-1,1}{}^{s_{n-1,n-\eps}}( \res_{n-\eps-1,1}^{n-\eps} N\boxtimes L^A(a^{(\eps)}))]\\
  &+[ \ind_{n-\eps,\eps-1,1}^{n-1,1}{}^{s_{n-1,0,n-\eps}}(N\boxtimes \res_{\H^A_{1,\eps-1}}^{\H^A_\eps}L^A(a^{(\eps)}))]
\end{split}\end{equation*}
The third subquotient does not contribute to $\De_a \ind_{n-\eps,\eps}^{n} N\boxtimes L^A(a^{(\eps)})$ as its formal character ends on $a^{-1}$. Similarly, the direct summands of $\res_{\H^A_{\eps-1,1}}^{\H^A_\eps} L^A(a^{(\eps)})$ other than $\De_a L^A(a^{(\eps)})$ and the direct summands of $\res_{n-\eps-1,1}^{n-\eps} N\boxtimes L^A(a^{(\eps)})$ other than $\De_a N\boxtimes L^A(a^{(\eps)})$ cannot contribute to $\De_a\ind_{n-\eps,\eps}^{n} N\boxtimes L^A(a^{(\eps)})$. 
As $\De_a N=0$, we obtain
$$
\De_a\ind_{n-\eps,\eps}^{n} N\boxtimes L^A(a^{(\eps)}) \cong  
\ind_{n-\eps,\eps-1,1}^{n-1,1} N\boxtimes \De_a L^A(a^{(\eps)}).
$$
By considering characters, we see that 
$$[\De_a L^A(a^{(\eps)})] = \eps [L^A(a^{(\eps - 1)}) \boxtimes (a)],$$
hence
\begin{equation*}\label{E300900_3}
[\De_a\ind_{n-\eps,\eps}^{n} N\boxtimes L^A(a^{(\eps)})] =  
\eps[\ind_{n-\eps,\eps-1,1}^{n-1,1} 
N\boxtimes L^A(a^{(\eps-1)})\boxtimes (a)].
\end{equation*}
By Lemma~\ref{sunny}(ii), 
the cosocle of $\ind_{n-\eps, \eps-1,1}^{n-1,1} N \boxtimes L^A(a^{(\eps-1)})
\boxtimes (a)$ is  $(\tilde f_a^{(\eps-1)} N) \boxtimes (a)$ which is the same 
as $(\tilde e_a M) \boxtimes (a)$, and all other composition
factors of this module are of the form $L \boxtimes  (a)$
with $\eps_a(L) < \eps - 1$ by Lemma~\ref{L010900}.
Moreover, all composition factors of $\De_a K$
are of the form $L \boxtimes  (a)$
with $\eps_a(L) < \eps - 1$.
So we have now seen that
$$
[\De_a M] = \eps [\tilde e_a M \boxtimes (a)] + \sum c_r [N_r 
\boxtimes  (a)]
$$
for irreducibles $N_r$ with $\eps_a(N_r) < \eps_a(\tilde e_a M)$, which 
implies (i).

The proofs of (ii) and (iii) are analogous to \cite{K:05}, Theorem 5.5.1.
\end{proof}

An interesting consequence is the following.

\begin{mycor}\label{epr}
Let $M, N\in\Rep_{\neq \pm 1} \H_n$ be irreducible modules with \mbox{$M \not\cong N$.}
Then, for $a\in F \setminus \{\pm 1\}$, we have 
$\hom_{\H_{n-1}}(e_a M, e_a N) = 0.$
\end{mycor}

\begin{proof}Analogous to \cite{BK:01}, Corollary 6.12 or \cite{K:05}, Corollary 5.5.2.
\end{proof}

\section{Discussion of eigenvalues $1$ and $-1$}

We now investigate the cases when some of the $X_i$ have
eigenvalues $1$ and $-1$. 
The usual formal character argument used in section $4$ does not fully apply here, since, when we compute formal characters of induced modules $f_a N$, the term coming from the longest coset respresentative in the Shuffle Lemma, ends on $a^{-1}$ which for $a = \pm 1$ is the same as $a$. Thus we obtain higher multiplicities of $N \boxtimes (a)$ in $\De_a f_a N$.
First, we consider eigenspaces with a maximal number of $a$s at the end. 

\begin{lemma}\label{eigspaces}
Let $N \in \H_{n-1} \mod$ be irreducible and $ a = \pm 1$.
Let $\eps_{a} (N) =r-1 \geq 0$. Set $M:= \ind_{n-1,1}^{n} N \boxtimes (a)$. 
Let $(a_0,\underline{b},a^{(r-1)})$ be a tuple of eigenvalues on $N$. Then any
$(a_0,\underline{b},a^{(r)})$-eigenvector in $M$ is either of the form $$1
\otimes v$$ or 
$$T_{n-r,0,n-1} \otimes v  +\sum_{j=0}^{n-r-1} T_{j,0,n-1}\otimes
   u_{j} + \sum_{j=1}^{n-1} T_{j,n-1} \otimes w_{j} + 1\otimes u_0 $$
where $v =\tilde v \boxtimes c$
for an  $(a_0,\underline{b},a^{(r-1)})$- \ or
   $(a_0a,\underline{b},a^{(r-1)})$- eigenvector $\tilde v$ in $N$, a
   generator $c$ of the one-dimensional module $(a)$ and some
$u_j, w_j,u_0$ in $N \boxtimes (a)$.

In other words, it has leading term $1$ or $T_{n-r,0,n-1} $.

In particular, the  $(a_0,\underline{b},a^{(r)})$-eigenspace in M has at
most twice the dimension of that in $N \boxtimes (a)$.
\end{lemma} 

\begin{proof}
Directly from the action of the Weyl group on the lattice, we see that the
leading term has to be of the form $T_{l,0,n-1}$ for $l \geq n-r$ or $T_{l,n-1}$
for $l \geq n-r+1$ for there to be an $a$-tail of length $r$. 

For leading term $T_{l,0,n-1}$ with $l > n-r$
\begin{equation*}\begin{split}(X_{l}-a)&(T_{l,0,n-1} \otimes u_l + T_{l-1,0,n-1} \otimes u_{l-1})\\
= &T_{l,0,n-1}  (X_{l}-a) \otimes u_{l}  - (q-q^{-1})T_{l-1,0,n-1} X_{l}
\otimes u_{l} \\ &+ T_{l-1,0,n-1} (X_n^{-1}-a) \otimes u_{l-1}+ \hbox{lower terms}
\end{split}\end{equation*}
Since $X_n^{-1}-a$ acts as $0$ on the whole of $N \boxtimes (a)$ and since none of the lower
terms can contribute to the coefficient of $T_{l-1,0,n-1}$, this shows that $u_l =0$
which contradicts the assumption that $T_{l,0,n-1}$ is the leading term.
Thus the leading term can only be $T_{n-r,0,n-1}$ in this case.

Equally, if the leading term is $T_{l,n-1}$ for $l \geq  n-r+1$ 
\begin{equation*}\begin{split}
(X_{l+1}-a)&(T_{l,n-1} \otimes w_l + T_{l+1,n-1} \otimes w_{l+1})\\
=&T_{l,n-1}  (X_{l}- a) \otimes w_{l}  + (q-q^{-1})T_{l+1,n-1} X_{n}
\otimes w_{l}\\ &+ T_{l+1,n-1} (X_n-a) \otimes w_{l+1}+ \text{lower terms}
\end{split}\end{equation*}
By the same argument as above, $w_l=0$ and the leading term has to be $1$.
\end{proof}

\begin{theorem}\label{endring}

Let $N \in \H_{n-1} \mod$ be irreducible, $a = \pm 1$ and set $$M:= f_aN
.$$ Then $ \dim \: \End_{\H_n}(M) \leq 2$. 

\end{theorem} 

\begin{proof}
By Frobenius reciprocity, $\End_{\H_n}(M) \cong \hom_{\H_{n-1,1}}( N \boxtimes (a),\De_a M)$. Let $r-1=\eps_{a} (N)$.
By Lemma \ref{eigspaces}, the  $(a_0,\underline{b},a^{(r)})$-eigenspace in M has at
most twice the dimension of that in $N \boxtimes (a)$, therefore the socle of
$\De_a M$ can contain at most 2 copies of $N \boxtimes (a)$.
\end{proof}

Recall the automorphism $\si$ fixing all generators except $X_0$ which it maps to its negative.

\begin{mycor}\label{eitheror}
If $N$ in Lemma \ref{endring} satisfies $N \cong N^\tau$, $M= 
f_a N$ either has an irreducible cosocle or splits into a direct sum
of two non-isomorphic irreducibles.
\end{mycor}

\begin{proof}
There are two cases to consider. The first case is the case $a=-1$. In this case, by $N \cong N^\tau$ and Corollary \ref{c2}, we have
\begin{equation*}\begin{split}
(\ind_{n-1,1}^{n} N \boxtimes (-1))^\tau &\cong \ind_{n-1,1}^{n}{^{s_{n-1,0,n-1}}( N^\tau \boxtimes (-1)^\tau)} \\ 
&\cong \ind_{n-1,1}^{n}{^{s_{n-1,0,n-1}} (N \boxtimes (-1))} \\ 
&\cong \ind_{n-1,1}^{n} N^\si \boxtimes (-1).
\end{split}\end{equation*}
(where $N^\si$ might or might not bei isomorphic to $N$).
Now $\tau$ fixes the $X_i$, duality doesn't affect eigenvectors and
$\soc f_{-1}N \cong (\tilde f_{-1}N^\si)^\tau$. Therefore the usual
argument of $\hom_{\H_{n-1,1}}(N^\si \boxtimes (-1), \De_{-1}Q)$ being
nonzero for every constituent of $\tilde f_{-1}N^\si$ together with the
restriction on the dimension of $(\pm
a_0,\underline{b},-1^{(\eps_{-1}(N))})$-eigenspaces from Lemma
\ref{eigspaces} yields that the socle of $f_{-1}N$ can have at most two constituents. In case it is irreducible, $\cosoc f_{-1}N^\si$ and therefore also $\cosoc f_{-1}N$ is irreducible.
In case $\soc f_{-1}N$ has two constituents, its $(\pm a_0,\underline{b},-1^{(\eps_{-1}(N))})$-eigenspace has to contain $1 \otimes N \boxtimes (-1)$ and therefore socle and cosocle have to coincide. Thus $f_{-1}N \cong M_1 \oplus M_2$ and by the restriction on the dimension of the endomorphism ring, these two irreducibles have to be nonisomorphic.

If $a=1$, then
$N \cong N^\tau$ implies $M \cong M^\tau$ since
\begin{equation*}\begin{split}
M^\tau &\cong (\ind_{n-1,1}^{n} N\boxtimes (a))^\tau\\
&\cong \ind_{n-1,1}^{n}{^{s_{n-1,0,n-1}}( N^\tau
  \boxtimes (a)^\tau)} \\ 
&\cong \ind_{n-1,1}^{n}{^{s_{n-1,0,n-1}} (N \boxtimes (a))} \\ 
&\cong \ind_{n-1,1}^{n} N \boxtimes (a),
\end{split}\end{equation*}
 so in fact $M$ is self-dual. Now the same argument about every constituent
 of the socle contributing
 $(a_0,\underline{b},-1^{(\eps_{-1}(N))})$-eigenspace as in the previous
 case shows that either the socle (and hence by self-duality the cosocle)
 must be irreducible or the socle must equal the cosocle. In the latter
 case the module itself is the direct sum of two non-isomorphic irreducibles.
\end{proof}

Now we want to consider a class of modules where
the results are as in the case $a \notin \{ \pm 1\}$.

\begin{lemma}\label{indsi}
For all $N \in \H_{n-1} \mod$ and all $a \in F$, 
 $$(\ind_{n-1,1}^{n} N \boxtimes (a))^\si \cong \ind_{n-1,1}^{n} N^\si \boxtimes (a).$$
\end{lemma}

\begin{proof}
Since $1 \otimes N\boxtimes (a)$ generates $\ind_{n-1,1}^{n} N \boxtimes (a)$ as well as
$(\ind_{n-1,1}^{n} N \boxtimes (a))^\si$ as an $\H_{n}$-module it suffices
to check that the twisted action on $1 \otimes N\boxtimes (a)$ is the same as
the action on $1 \otimes N^\si\boxtimes (a)$.
So let $h \in \H_n$ and write $h= \sum_{w \in D_{n-1,1}} T_w h_w$ for $h_w \in \H_{n-1,1}$ .
Then, as $\si$ fixes all of $\H_n^\fn$, 
\begin{equation*}\begin{split}
\si (h) \otimes u &= \si (\sum_{w \in D_{n-1,1}} T_w h_w)\otimes u\\
&=\sum_{w \in D_{n-1,1}} T_w \si (h_w)\otimes u\\
&=\sum_{w \in D_{n-1,1}} T_w \otimes \si (h_w) u
\end{split}\end{equation*}
for $u \in N \boxtimes (a)$, which proves the claim.
\end{proof}

\begin{theorem}\label{irred-1}
  Let $M \in \H_{n}\mod$ be irreducible with \mbox{$\ch M \neq \ch
    M^\si$} and let \mbox{$\eps_{-1}(M) = r \geq 1$.} Assume that
  there exists an irreducible submodule $N \boxtimes (-1)$ of
  $\De_{-1}M$ such that $N^\si \ncong N = L(a_0, \dots, a_{n-r},-1^{(r-1)})$ for some $(a_0, a_1, \dots, a_{n-r}) \in (F^\times \setminus \{\pm1\}
)^{n-r+1}$.
  Then
\begin{enumerate}
\item $\De_{-1^{(r)}} M \cong K \boxtimes L^A(-1^{(r)})$ for some irreducible
  $K\in \H_{n-r} \mod$ with \mbox{$\ch K \neq \ch K^\si$};
\item $M$ is the irreducible cosocle of $\ind_{n-r,r}^n K \boxtimes L^A(-1^{(r)})$;
\item $\soc \De_{-1} M \cong  N\boxtimes (-1)$, i.e. $\tilde e_{-1}M \cong N$;
\item $M$ is the irreducible cosocle of the non-irreducible module $f_{-1}N$, in particular $\tilde f_{-1}N \cong M$;
\item $M$ is uniquely determined by $\ch M$ and
  in particular \mbox{$M \cong M^\tau$.}
\end{enumerate}
\end{theorem}

\begin{proof}
Proceed by induction on $r$. The case $r=1$ is easy: Take an irreducible submodule $N \boxtimes (-1)$ of $\De_{-1}M$ with $N^\si \ncong N \cong L(a_0, \dots, a_{n-1})$. Then this will also play the role of $K$ in the statement and by the Mackey Theorem $$[\De_{-1}f_{-1} N] = [N \boxtimes (-1)]+[N^\si \boxtimes (-1)],$$ whence $\tilde f_{-1}N$ is irreducible. On the other hand $f_{-1} N$ is not irreducible since $\ch M \neq \ch M^\si$ and  $\De_{-1} M \cong N \boxtimes (-1)$ since some composition factor (and thus the composition factor $N^\si \boxtimes (-1)$) has to be taken up by $\soc f_{-1}N$ by the same argument as in Lemma \ref{eitheror}. $N$ is uniquely determined by $\ch N$ as $N \in \Rep_{\neq \pm 1} \H_{n-1}$, and if $\ch M = \ch M'$, $\De_{-1} M' \cong N \boxtimes (-1)$, whence $M' \cong \tilde f_{-1}N \cong M$.

Now inductively assume that (i) -(v) hold for $N$ in the hypothesis of the
theorem, in particular that there exists some irreducible   \mbox{$K \ncong K^\si \in
  \H_{n-r} \mod$} determined by $\ch K$ such that $\De_{-1^{(r-1)}} N \cong K \boxtimes
L^A(-1^{(r-1)})$, $N \cong \ind_{n-r,r-1}^{n-1} K \boxtimes L^A(-1^{(r-1)})$ and $\ch N$ determines $N$.

For the inductive step, we apply the Shuffle Lemma  to compute a filtration of \mbox{$\De_{-1^{(r)}} \ind_{n-1,1}^n N\boxtimes (-1)$.}
\mbox{$D_{(n-r,r),(n-1,1)} = \{1, s_{n-r,n-1}, s_{n-r,0,n-1} \}$}, but
considering formal characters we see that the subquotient of
$\res^n_{n-r,r} \ind_{n-1,1}^n N\boxtimes (-1)$ corresponding to
$s_{n-r,n-1}$ does not have formal character values with an $(-1)$-tail of
length $r$, as it is isomorphic to \mbox{$\ind_{n-r-1,1,r}^{n-r,r}
  {}^{s_{n-r,n-1}}(\res^{n-1,1}_{n-r-1,r,1}N\boxtimes (-1))$.} Thus it does not contribute to $\De_{-1^{(r)}} \ind_{n-1,1}^n N\boxtimes (-1)$ and the latter only has subquotients isomorphic to 
\begin{equation}\label{first}\ind_{n-r,r-1,1}^{n-r,r}(\De_{-1^{(r-1)}}N )\boxtimes (-1)\end{equation} 
and
\begin{equation}\label{second}\ind_{n-r,1,r-1}^{n-r,r} {}^{s_{n-r,0,n-1}}(\De_{-1^{(r-1)}}N)\boxtimes (-1).\end{equation}
The use of $\De_{-1^{(r-1)}}$ instead of $\res^{n-1,1}_{n-r,r-1,1}$ in
the Mackey formula is validated by the fact that no other summand of
$\res^{n-1}_{n-r,r-1}N $ can contribute composition factors to
\mbox{$\De_{-1^{(r)}} \ind_{n-1,1}^n N\boxtimes (-1)$} for lack of
$(-1)$s at the end.  Now, by the assumptions on $N$, (\ref{first}) is
isomorphic to $ K \boxtimes L^A(-1^{(r)})$ while (\ref{second}) has the
same formal character as $ K^\si \boxtimes L^A(-1^{(r)})$ and is
therefore isomorphic to $ K^\si \boxtimes L^A(-1^{(r)})$.  As in Lemma \ref{eitheror}, the
assumption on $N$ and Corollary \ref{c2} yield $$(\ind_{n-1,1}^{n} N \boxtimes (-1))^\tau \cong \ind_{n-1,1}^{n} N^\si \boxtimes (-1).$$
Now, if $\ind_{n-1,1}^n N\boxtimes (-1)$ were irreducible and therefore isomorphic to $M$, $M^\tau$ would also have to be isomorphic to $\ind_{n-1,1}^n N^\si \boxtimes (-1)$ and therefore $M^\tau \cong M^\si$, a contradiction since $$\ch M^\tau = \ch M \neq \ch M^\si.$$ Therefore, $\ind_{n-1,1}^n N\boxtimes (-1)$ is not irreducible.
But since, for any quotient $Q$ of $\ind_{n-1,1}^n N\boxtimes (-1)$ and in
particular for any constituent of \mbox{$\cosoc \ind_{n-1,1}^n N\boxtimes
(-1)$}, 
\begin{equation*}\begin{split} \hom_{\H_n}(\ind_{n-r,r}^n K \boxtimes L^A(-1^{(r)}) &, Q) \\&\cong
\hom_{\H_{n-r,r}}(K \boxtimes L^A(-1^{(r)}) ,\De_{-1^{(r)}}  Q)\end{split}\end{equation*} is nonzero and at the same time at most one-dimensional (as $K \boxtimes L^A(-1^{(r)})$ appears only once in $\De_{-1^{(r)}} \ind_{n-1,1}^n N\boxtimes
(-1)$), this shows that 
$\cosoc \ind_{n-1,1}^n N\boxtimes (-1)$ is irreducible and therefore isomorphic to
$M$, which proves (iv).

The same argument applied to $\ind_{n-1,1}^{n} N^\si \boxtimes (-1)$ shows
that $$\widetilde M = \soc \ind_{n-1,1}^n N\boxtimes (-1) \cong \cosoc \ind_{n-1,1}^n N^\si \boxtimes (-1)$$ is also irreducible
and has $\eps_{-1} (\widetilde M)=r$. Therefore $\widetilde M$ contributes the
composition factor $ K^\si \boxtimes L^A(-1^{(r)})$ to $\De_{-1^{(r)}}
\ind_{n-1,1}^n N\boxtimes (-1)$ and $$\De_{-1^{(r)}} M \cong K \boxtimes
L^A(-1^{(r)}),$$ so (i) holds.

To see that $N\boxtimes (-1)$ is actually the whole of $\soc \De_{-1} M$ it
suffices to consider the fact that, if $N'\boxtimes (-1)$ is an irreducible
submodule of $\De_{-1} M$, it has to satisfy $\eps_{-1}(N') =r-1$ as in this
case $$\hom_{\H_n}(\ind_{n-1,1}^{n} N'\boxtimes (-1),M) \neq 0.$$ But then $0
\neq \De_{-1^{(r-1)}}N' \cong K'\boxtimes L^A(-1^{(r-1)})$ whence, by
transitivity of induction,  $\ind_{n-r,r}^n K'\boxtimes L^A(-1^{(r)})$
projects onto $M$ or equivalently $$K' \boxtimes L^A(-1^{(r)}) \hookrightarrow
\De_{-1^{(r)}} M,$$ a contradiction, whence (iii) follows.

Now, suppose some $M'$ has the same formal character as $M$. Then $\De_{-1^{(r)}} M'$ has
the same formal character as $\De_{-1^{(r)}} M$. By virtue of $K$ being uniquely determined by its formal character, we obtain
$$\De_{-1^{(r)}} M' \cong K \boxtimes L^A(-1^{(r)}),$$ so (v) will follow if we
can show (ii), namely that $\cosoc \ind_{n-r,r}^n K \boxtimes L^A(-1^{(r)})$
is irreducible. 
It certainly does not contain several copies of $M$ since 
\begin{equation*}\begin{split}
\hom_{\H_{n-r,r}}(  K \boxtimes L^A(-1^{(r)}) &,\De_{-1^{(r)}}  M)\\& \cong
 \hom_{\H_n}(\ind_{n-r,r}^n K \boxtimes L^A(-1^{(r)}) , M)
\end{split}\end{equation*}
 is one-dimensional, so assume it contains an irreducible submodule $M'$ not
 isomorphic to $M$. But again, choosing $N' \boxtimes (-1)$ in the socle of $\De_{-1} M'$ shows that necessarily $$\De_{-1^{(r-1)}} N' \cong K \boxtimes L^A(-1^{(r-1)}),$$ whence (by the irreducibility of $\cosoc \ind_{n-r,r-1}^{n-1} K \boxtimes L^A(-1^{(r-1)})$)  $N' \cong N$ and we're done.
\end{proof}

\begin{remark}
By an additional induction over the number of $-1$-strings in the tuple labelling $N$ in Theorem \ref{irred-1}, which mainly uses that applying $\tilde f_a$ for $a \neq \pm1$ in between the $-1$-strings won't affect the uniqueness of formal characters, we can omit the hypothesis that the $a_i$ ($1 \leq i \leq n-r$) are all distinct from $-1$. 
\end{remark}

\begin{mycor}
In the situation of Theorem \ref{irred-1}, $$\soc \ind_{n-1,1}^n N\boxtimes (-1) \cong \cosoc \ind_{n-1,1}^n N^\si \boxtimes (-1) \cong M^\si.$$
\end{mycor}



We have seen that for $a=-1$, as long as $M \ncong M^\si$, everything works as in the case $a \neq \pm 1$, and $M \ncong M^\si$ is satisfied until, for the first time, the induced module $\ind_{n-1,1}^n N \boxtimes (-1)$ is irreducible and therefore isomorphic to its $\si$-conjugate. 
The next lemma explicitly computes eigenspaces in induced modules.

\begin{lemma}
Let  $N \cong N^\tau \in \H_{n-1} \mod$ be irreducible, $a \in \{ \pm 1\}$ and
$\eps_{a} (N) =r-1$. Set $M:=f_a N
$.
The generalized $(a^{(r)})$-eigenspace for $X_{n-r+1}, \dots, X_n$ in $M$ is contained in the socle and
cosocle of $M$, thus all other composition factors $K$ have $\eps_{a} (K) \leq r-1$ .
\end{lemma} 

\begin{proof}
By the Mackey Theorem \ref{TMackey} ,
\begin{equation*}\begin{split}[\De_{a^{(r)}}M] &= [\ind_{n-r,r-1,1}^{n-r,r}
\De_{a^{(r-1)}}N \boxtimes (a)  ]\\ & \quad+[\ind_{n-r,1,r-1}^{n-r,r}
{^{s_{n-r,0,n-1}}(\De_{a^{(r-1)}}N \boxtimes (a))}]\end{split}\end{equation*}
and as in Lemma \ref{eigspaces} we know that the $(a^{(r)})$-eigenspace of \mbox{$(X_{n-r+1}, \dots , X_n)$} in the submodule $\ind_{n-r,r-1,1}^{n-r,r}
\De_{a^{(r-1)}}N \boxtimes (a)$, which is by Theorem \ref{TMackey} contained
in the socle of $\De_{a^{(r)}}M$, is contained in $1
\otimes \De_{a^{(r-1)}}N \boxtimes (a)$. Under the restriction of the
projection of $M$ onto its cosocle, this certainly doesn't map to zero,
whence we have an injection of $\ind_{n-r,r-1,1}^{n-r,r}
\De_{a^{(r-1)}}N \boxtimes (a)$ into $\De_{a^{(r)}}\cosoc M$.
Since  $M$ is self-dual in case $a=1$ or $a=-1$ and $N \cong N^\si$ and $M^\tau \cong M^\si$ and therefore $\soc M \cong (\cosoc M^\si)^\tau$ in the remaining case, the generalized $(a^{(r)})$-eigenspace $\De_{a^{(r)}}\soc M$ of the socle must have the same dimension, thus exhausting all of $\De_{a^{(r)}}M$.
\end{proof}

We will now give an example where indeed the functor $\tilde f_1$ does not
produce an irreducible module. 
\begin{example}\label{smallcounterex}
We apply the functor $f_1$ to the two-dimensional irreducible module
$L(a_0,q^2)\cong\ind_{\P_1}^{\H_1}(a_0,q^2) \in \H_1 \mod$ with basis $\{w_1,w_2\}$ on which $T_0, X_0,X_1$ act by matrices 
$$\begin{array}{ccc}
\begin{pmatrix} 0&1\\1&(p-p^{-1})
\end{pmatrix}, &
\begin{pmatrix} a_0&-(p-p^{-1})a_0q^2\\0&a_0q^2
\end{pmatrix}, &
\begin{pmatrix} q^2 & (p-p^{-1})(q^2+1)\\0&q^{-2}
\end{pmatrix}
\end{array}$$
respectively. 
Since $w_1 \in L(a_0,q^2)$ is an $(a_0,q^2)$-eigenvector, Lemma
\ref{intertwining} implies that, in the induced module 
$M:=f_1 L(a_0,q^2)  = \H_2
\otimes_{\H_{1,1}}  L(a_0,q^2) \boxtimes (1)$, we find an $(a_0,1,q^2)$-eigenvector 
$v_1:=(T_1 + q^{-1})w_1$. Again by Lemma \ref{intertwining}, the vector \mbox{$
  (T_0 - q^2T_0^{-1})w_1 \in  L(a_0,q^2)$}
 is an
$(a_0q^2,q^{-2})$-eigenvector for $(X_0,X_1)$, and \mbox{$v_2:=(T_1-q)(T_0 - q^2T_0^{-1})w_1$} is an
$(a_0q^2,1,q^{-2})$-eigenvector in M. 
It is then easy to calculate that $v_1$ and $v_2$ both generate $4$-dimensional
submodules with trivial intersection which are both irreducible and isomorphic
to $L(a_0,1,q^2)$ and $L(a_0q^2,1,q^{-2})$ respectively. Therefore $$\tilde f_1
L(a_0,q^2) \cong f_1
L(a_0,q^2) \cong L(a_0,1,q^2) \oplus L(a_0q^2,1,q^{-2}).$$
\end{example}

\section{Facts About Representations of the Affine Hecke Algebra of Type $A$}

At this point it is convenient to recall some important facts from the
representation theory of affine Hecke algebras of type $A$. The results in this
section are compiled from \cite{BZ:77}, \cite{A:96}, \cite{V:02}, \cite{V:99}, \cite{GV:01} and  \cite{G:99}.

For nonzero $\la \in F$, define the full subcategory $\Rep_\la \H_n^A$ of $\H_n^A \mod$ 
to consist of those modules where all eigenvalues of the generators $\R_n$ are from
$I_\la^+ := \{ \la q^{2i} | i \in \ZZ    \}$. The significance of the "plus" will
become clear later on when we move to type $B$.

The $\Rep_\la \H_n^A$ are equivalent categories for all $\la \in
F$, and for $I_\la \neq I_{\la^\prime}$, there are no nontrivial extensions between modules in $\Rep_\la \H_n^A$
and $\Rep_{\la^\prime} \H_n^A$ for $I_\la \neq I_{\la^\prime}$ and $\ind_{ \H_{n_1}^A
\otimes \H_{n_2}^A}^{ \H_{n_1+n_2}^A} M \boxtimes N$, for irreducible $M$ and
$N$ in $\Rep_\la \H_{n_1}^A$ and $\Rep_{\la^\prime} \H_{n_2}^A$ respectively, is
always irreducible.
The irreducible modules in $\Rep_\la \H_n^A$ are well-understood and have a
nice combinatorial description.

Call the sequence $\Ga_{(i, i+k)}=(aq^{2i}, aq^{2(i+1)}, \dots, aq^{2(i+k)})$ a segment and denote by $L_{\Ga_{(i, i+k)}} = L^A(aq^{2i}, aq^{2(i+1)}, \dots, aq^{2(i+k)})$ the 
one-dimensional  representation of $\H^A_{k+1}$ on which all $T_j$, for $1 \leq j
\leq k$ act as $q$ and $X_l$ acts as $aq^{2(i+l-1)}$ for $1 \leq j
\leq k+1$. Define a multisegment  $\Ga=(\Ga_1, \dots, \Ga_m)$
to be a concatenation of several segments and denote by $L_\Ga:= L_{\Ga_1} \boxtimes \cdots \boxtimes L_{\Ga_m}$ the
one-dimensional representation for the tensor product of the corresponding
algebras. The length of a multisegment is the sum of the lengths of the
contained segments.
We have two different
orderings on multisegments, the so-called right and left orders.

In the right order, $\Ga_{(i,i+k)} > \Ga_{(j,j+l)}$ if
$i>j$ or  if $i=j$ and $l>k$.

In the left order, $\Ga_{(i,i+k)} \succ \Ga_{(j,j+l)}$ if $i+k > j+l$ or $i+k
= j+l$ and $j>i$.


Bernstein and Zelevinski \cite{BZ:77} showed that there is a one-to-one
correspondence between ordered multisegments of length $n$ and irreducible modules for the
affine Hecke algebra $\H_n^A$ of type $A$. This correspondence is given by inducing $L_\Ga$ for a  multisegment $\Ga$ up to
$\H^A_n$ and taking the -- always irreducible -- cosocle of this induced
module which is independent of whether we have chosen $\Ga$ in right or left order. 
Denote this cosocle by $N_\Ga$. 

There is also a combinatorial description of some form of branching
rules. Since the index of $\H_{n-1}^A$ in $\H_n^A$ is infinite, we substitute
normal induction by functors 
$$\ind_a^A := \ind_{\H_{n-1,1}^A}^{ \H_{n}^A} - \boxtimes (a): \Rep_\la \H_{n-1}^A
\longrightarrow  \Rep_\la \H_{n}^A $$ for every $a \in I_\la^+$. 
For irreducible $N$ in $\Rep_\la \H_{n-1}^A$, the cosocle of $\ind_a^A N$ is
irreducible and we denote this by $\tilde{f}_a^A N$.
Dually, we define $$\tilde{e}_a^A N := \res^{ \H_{n-1,1}^A}_{\H_{n-1}^A} \soc
\De_a N$$ where $\De_a N$ is the generalized $a$-eigenspace of $X_n$ in $N$ as
in type $B$. 
For irreducible $N$, $\tilde{e}_a^A N$ is again an irreducible module.
As in type $B$, we define $\eps_a(N)$ to be the largest number $r$ such that 
$\De_{a^{(r)}} N \neq 0$. 
Since the definitions of $\De_a$ and $\eps_a$ involve the action of the
lattice, we keep the notation from type $B$ and do not mark the symbols with an $^A$.
Analogously, define $$\ind_a^{*A} := \ind_{\H_{1,n-1}^A}^{ \H_{n}^A} (a) \boxtimes -: \Rep_\la \H_{n-1}^A
\rightarrow  \Rep_\la \H_{n}^A, $$ $$\tilde{f}_a^{*A} N := \cosoc
\ind_a^{*A} N \quad \hbox{ and } \quad \tilde{e}_a^{*A} N := \res^{ \H_{1,n-1}^A}_{\H_{n-1}^A} \soc
\De^*_a N,$$ where $\De_a^* N$ is the generalized $a$-eigenspace of $X_1$ on
$N$ which is, just as $\De_a N$, an $\H_{1,n-1}^A$-submodule of $N$
since $X_1$ commutes with $T_2, \dots, T_n$. Both $\tilde{f}_a^{*A} N$ and
$\tilde{e}_a^{*A} N$ are irreducible if $N$ is.
Lastly, $\eps^*_a(N)$ is defined as the maximal $r$ such that $\De^*_{a^{(r)}} N \neq 0$.

For the irreducible module $N_\Ga$ in $\Rep_\la \H_{n}^A $ there are
combinatorial algorithms to compute $\tilde{f}_a^A N_\Ga,\tilde{f}_a^{*A} N_\Ga
,\tilde{e}_a^A N_\Ga,\tilde{e}_a^{*A} N_\Ga , \eps_a(N_\Ga)$ and
$\eps^*_a(N_\Ga)$ which will be described below.

To compute $\eps_a(N_\Ga)$, write down the multisegment $\Ga$ in right
order. Then, write a $+$ for every segment ending on $aq^{-2}$ and a $-$ for
every multisegment ending on $a$. In the resulting sequence of plus and minus
signs successively cancel out all subsequences of the form $-+$ until the
leftover sequence is of the form $+ \cdots + -\cdots -$.
The number of uncanceled $-$ signs is $\eps_a(N_\Ga)$.
If we replace the segment $(aq^{-2i}, \dots,aq^{-2}  , a)$ which contributed the
leftmost uncanceled $-$, by $(aq^{-2i}, \dots, aq^{-2})$ we get the
multisegment corresponding to $\tilde{e}_a^A N_\Ga$. Denote this multisegment
by $\tilde{e}_a^A \Ga$. In case there is no $-$ sign left after cancellation
\mbox{$\tilde{e}_a^A N_\Ga = 0$.}
If we replace the segment $(aq^{-2k}, \dots,aq^{-2} )$ which contributed the
rightmost uncanceled $+$, by $(aq^{-2k}, \dots, aq^{-2},a)$ we get the
multisegment $\tilde{f}_a^A \Ga$ corresponding to $\tilde{f}_a^A N_\Ga$.
If there is no $+$ left after cancellation, we add a new segment $(a)$ to $\Ga$
to obtain $\tilde{f}_a^A N_\Ga$.

To compute $\eps^*_a(N_\Ga)$, write down the multisegment $\Ga$ in left
order. Then, write a $+$ for every segment starting on $aq^{2}$ and a $-$ for
every multisegment starting on $a$. In the resulting sequence of plus and minus
signs successively cancel out all subsequences of the form $+-$ until the
leftover sequence is of the form $- \cdots - +\cdots +$.
The number of uncanceled $-$ signs is $\eps^*_a(N_\Ga)$.
If we replace the segment $(a,aq^2,\dots  , aq^{2i})$ which contributed the
leftmost uncanceled $-$, by $(aq^{2}, \dots, aq^{2i})$ we get the
multisegment corresponding to $\tilde{e}_a^{*A} N_\Ga$. Denote this multisegment
by $\tilde{e}_a^{*A} \Ga$. If there is no $-$ sign left after cancellation 
\mbox{$\tilde{e}_a^{*A} N_\Ga = 0$.}
If we replace the segment $(aq^2, \dots,aq^{2k} )$ which contributed the
rightmost uncanceled $+$, by $(a,aq^2, \dots, aq^{2k})$ we get the
multisegment $\tilde{f}_a^{*A} \Ga$ corresponding to $\tilde{f}_a^{*A} N_\Ga$.
If there is no $+$ left after cancellation, we add a new segment $(a)$ to $\Ga$
to obtain $\tilde{f}_a^{*A} N_\Ga$.
Moreover, $$\tilde{f}_a^{*A} N_\Ga \cong \soc \ind_a^A N_\Ga \quad \hbox{ and } \quad \tilde{f}_a^{A} N_\Ga \cong \soc \ind_a^{*A} N_\Ga.$$

\section{Subcategories with Type-$A$ Behavior}\label{easy}

In this Chapter we inverstigate certain subcategories of $\H_n \mod$ which behave very similarly to the situation in type $A$.
We define $\Rep_\la \H_n$ for fixed nonzero $\la \in F$ to be the full subcategory of
$\H_n^R \mod$ where all
eigenvalues of $\R_n$ are from the set $$I_\la := \{ \la q^{2i},\la^{-1}q^{2i}
| i \in \ZZ  \} .$$

In this section, we consider the cases where $p^2, \pm q, \pm1 \notin I_\la$.
In these cases 
$I_\la^+ := \{\la q^{2i} | i \in \ZZ \}$ and $I_\la^- := \{\la^{-1} q^{2i} | i \in \ZZ \}$
are disjoint.
Since we exclude the eigenvalue $-1$, Lemmas \ref{cliff} and  \ref{cliff2}
guarantee that we can work with the algebra $\H_n^R$ instead of $\H_n$ since all irreducibles for $\H_n$ in the analogous subcategories are obtained by extending the action to $\H_n$ with an arbitrary new eigenvalue for $X_0$. Working with $\H_n^R$ here is more
convenient since we want to exploit the subalgebra $\H_n^A$, which has finite
index in $\H_n^R$ but not in $\H_n$.

First, we would like to give a very general result on the
formal characters of an $\H_n^R$-module obtained by inducing from $\H_n^A$.

\begin{lemma}\label{indAB}
Let $N \in \H_n^A \mod$ be irreducible. Set $M:= \ind_{
  \H_n^A}^{\H_n^R}  N.$  Then $$\eps_a(M) \leq \eps_a(N) + \eps^*_{a^{-1}}(N).$$
\end{lemma}
\proof
By the Shuffle Lemma, we obtain the formal character of $M$ by taking all
formal characters of $N$ and then successively inverting the first
entries of the tuple and moving
them to the rear, see Section 3 for an example.
From this, one easily sees that the maximal number of $a$'s at the end of a
tuple from $\ch M$ is less
or equal to the maximal number of $a$'s at the end of a
tuple from $\ch N$ plus the maximal number of $a^{-1}$'s at the beginning of a
tuple from $\ch N$ that get inverted and moved to the rear.
\endproof

Note that on $\H_n^A$ there is an algebra antiautomorphism $\kappa$ given by 
$T_i \mapsto T_{n-i}$ and $X_i \mapsto X_{n+1-i}^{-1}$ inducing a duality on
$\H_n^A \mod$. It is easy to check on the generators that 
$\kappa$ is the composite of first
taking the $\tau$-dual and then twisting with the longest coset representative 
$d= s_0s_{1,0} \cdots s_{j,0} \cdots s_{n-1,0}$.
Now we compute the $\kappa$-dual of an irreducible $N_\Ga \in \Rep_\la
\H_n^A$, where $\Ga$ consists of segments $\Ga_1, \dots, \Ga_r$ of length
$n_1, \dots, n_r$ respectively:
\begin{equation*}\begin{split}
N_\Ga^\kappa &\cong {^d(N_\Ga^\tau)}\\
&\cong {^d(N_\Ga)} \\
&= \cosoc {^d(
 \H_n^A \otimes_{\H_{n_1}^A \otimes \cdots \otimes \H_{n_r}^A}L_\Ga)} \\
&\cong \cosoc \H_n^A \otimes_{^d(\H_{n_1}^A \otimes \cdots \otimes \H_{n_r}^A)}{^dL_\Ga}\\
&\cong \cosoc \H_n^A \otimes_{\H_{n_r}^A \otimes \cdots \otimes \H_{n_1}^A}
L_{\overline{ \Ga}}\\
&\cong N_{\overline{ \Ga}},
\end{split}\end{equation*}
where $\overline{ \Ga}$ is the multisegment $\overline{ \Ga}_r, \dots,
\overline{ \Ga}_1$, and the segment $\overline{ \Ga}_j$, for a segment $\Ga_j =(a, \dots, aq^{2k})$,
 is defined as $\overline{ \Ga}_j= (a^{-1}q^{-2k}, \dots , a^{-1})$.
The second isomorphism uses the fact that in type $A$ all irreducibles are self-dual under $\tau$-duality.

\begin{lemma}\label{easycase}
Let $N_\Ga \in \Rep_{\la^{-1}} \H_n^A$ be irreducible. Then $M_\Ga:=\ind_{\H_n^A}^{\H_n^R}
N_\Ga \in \Rep_{\la} \H_n$ 
is irreducible.
\end{lemma} 

\begin{proof}
Since all formal characters of $N_\Ga$ have entries in $I_\la^-$ and $I_\la^-
\cap I_\la^+ = \emptyset$, the Shuffle Lemma implies that the only summands in
the formal character of $M_\Ga$
 exclusively containing entries from 
$I_\la^-$ are those afforded by the coset representative $1$.
Thus $$\hom_{\H_n^R}(M_\Ga, \cosoc M_\Ga) \cong \hom_{  \H_n^A}( N_\Ga, \res^{\H_n^R}_{  \H_n^A} \cosoc M_\Ga ) \cong F,$$
whence the cosocle of $M_\Ga$ is irreducible.
Since by Corollary \ref{c2}
\begin{equation*}\begin{split}\hom_{\H_n^R}&(\soc M_\Ga, M_\Ga) \cong  \hom_{\H_n^R}(M_\Ga^\tau, \cosoc
M_\Ga^\tau) \\ 
& \cong \hom_{\H_n^R}(\ind_{  \H_n^A}^{\H_n^R} {^d(
  N_\Ga^\tau)}, \cosoc\ind_{  \H_n^A}^{\H_n^R} {^d(N_\Ga^\tau)} ) \\ 
&\cong \hom_{\H_n^R}(\ind_{  \H_n^A}^{\H_n^R} N_{\overline{\Ga}},\cosoc\ind_{  \H_n^A}^{\H_n^R} N_{\overline{\Ga}})\\
& \cong \hom_{  \H_n^A}(N_{\overline{\Ga}}, \res^{\H_n^R}_{  \H_n^A}  \cosoc\ind_{  \H_n^A}^{\H_n^R} N_{\overline{\Ga}}   )\\
& \cong F
\end{split}\end{equation*}
by the same argument as above, the socle of  $M_\Ga$ is also simple and
contains the generalized simultaneous eigenspace of $\R_n$ where all
eigenvalues of the $X_i$ are from $ I_\la^+$.

Now take a $(b_1, \dots, b_n)$-eigenvector $u \in M_\Ga$ such that $b_1, \dots, b_n$
are all in $ I_\la^+$. 
Then 
\begin{equation*}\begin{split}
v:=&(T_0 -b_nT_0^{-1})(T_{1} -b_{n-1}b_nT_{1}^{-1})(T_0 -b_{n-1}T_0^{-1}) \\
& \cdots (T_{n-j} -b_jb_nT_{n-j}^{-1})\cdots(T_0 -b_jT_0^{-1}) \\
&\cdots (T_{n-1} -b_1b_nT_{n-1}^{-1})\cdots(T_1 -b_1b_2T_1^{-1}) (T_0 -b_1T_0^{-1})u
\end{split}\end{equation*}
is a $(b_n^{-1}, \dots, b_1^{-1})$-eigenvector  by Lemma \ref{intertwining}
since $b_j \notin \{ p^{\pm 2},1\}$ and $b_jb_k \notin \{q^{\pm 2}\}.$
Now all
$b_i^{-1}$ belong to  $ I_\la^-$, therefore $v \in 1\otimes N_\Ga$ 
and generates $M_\Ga$. Therefore any element in the socle of $M_\Ga$
generates the whole of $M_\Ga$, so it must be irreducible.
 \end{proof}

\begin{lemma}Let $N_\Ga \in \Rep_{\la^{-1}} \H_n^A$ be irreducible. 

Then
$\eps_a(\ind_{  \H_n^A}^{\H_n^R} N_\Ga   ) = \left\{
\begin{array}{ll}
\eps_a(N_\Ga) & \hbox{if } a \in I_\la^- \\
\eps^*_{a^{-1}}(N_\Ga) & \hbox{if } a \in I_\la^+ 
\end{array} \right.$
\end{lemma} 

\begin{proof}
This follows directly from the proof of Lemma \ref{indAB} and the fact that  $I_\la^-
\cap I_\la^+ = \emptyset$.
\end{proof}

\begin{lemma}
Let $M_\Ga \in Rep_\la \H_n$ be defined as in Lemma \ref{easycase} and $a \in I_\la$. Then 
\begin{enumerate}
\item
$\tilde{f}_a M_\Ga =\left\{
\begin{array}{ll}
M_{\tilde{f}_a^A \Ga} & \hbox{if } a \in I_\la^- \\
M_{\tilde{f}_{a^{-1}}^{*A} \Ga} & \hbox{if } a \in I_\la^+
\end{array} \right.$

\item
$\tilde{e}_a M_\Ga =\left\{
\begin{array}{ll}
M_{\tilde{e}_a^A \Ga} & \hbox{if } a \in I_\la^- \\
M_{\tilde{e}_{a^{-1}}^{*A} \Ga} & \hbox{if } a \in I_\la^+
\end{array} \right.$

\end{enumerate}
\end{lemma} 

\begin{proof}
(i) Without loss of generality we assume $a \in I_\la^-$ and compute $\tilde{f}_a
M_\Ga$ and $\tilde{f}_{a^{-1}} M_\Ga = \soc \ind_{H_{n,1}^R}^{\H_{n+1}^R} M_\Ga \boxtimes (a)$.

We know that
\begin{equation*}\begin{split}
\ind_{H_{n,1}^R}^{\H_{n+1}^R} M_\Ga \boxtimes (a)
& \cong\ind_{H_{n,1}^R}^{\H_{n+1}^R} (\ind_{H_{n}^A}^{\H_{n}^R}  N_\Ga )\boxtimes (a)\\
& \cong \ind_{H_{n+1}^A}^{\H_{n+1}^R} \ind_{H_{n,1}^A}^{\H_{n+1}^A} N_\Ga\boxtimes (a),
\end{split}\end{equation*}
so, as every composition factor in $\ind_{H_{n,1}^A}^{\H_{n+1}^A} N_\Ga\boxtimes
(a)$ is in $\Rep_{\la^{-1}} \H_{n+1}^A$ and therefore yields an irreducible
subquotient of $\ind_{H_{n,1}^R}^{\H_{n+1}^R} M_\Ga \boxtimes (a)$ upon
induction to type $B$, $\ind_{H_{n,1}^R}^{\H_{n+1}^R} M_\Ga \boxtimes (a)$ has the
same number of composition factors as  $\ind_{H_{n,1}^A}^{\H_{n+1}^A} N_\Ga\boxtimes
(a)$, labeled by the same multisegments.
Since socle and cosocle of  $\ind_{H_{n,1}^R}^{\H_{n+1}^R} M_\Ga \boxtimes (a)$
are irreducible, they have to coincide with  $$\ind_{H_{n+1}^A}^{\H_{n+1}^R} \soc
\ind_{H_{n,1}^A}^{\H_{n+1}^A} N_\Ga\boxtimes (a)$$ and $$\ind_{H_{n+1}^A}^{\H_{n+1}^R} \cosoc
\ind_{H_{n,1}^A}^{\H_{n+1}^A} N_\Ga\boxtimes (a)$$ respectively. The fact that
$N_{\tilde{f}_{a}^{*A} \Ga} = \soc \ind_{H_{n,1}^A}^{\H_{n+1}^A} N_\Ga\boxtimes
(a)$ implies that $$\tilde{f}_{a^{-1}} M_\Ga \cong  \ind_{H_{n+1}^A}^{\H_{n+1}^R}N_{\tilde{f}_{a}^{*A} \Ga}$$ for $a \in I_\la^-$, which completes the proof of (i).

(ii) follows directly from Lemma \ref{L290900_2}.
\end{proof}

\bibliographystyle{amsplain}

\providecommand{\bysame}{\leavevmode\hbox to3em{\hrulefill}\thinspace}
\providecommand{\MR}{\relax\ifhmode\unskip\space\fi MR }
\providecommand{\MRhref}[2]{%
  \href{http://www.ams.org/mathscinet-getitem?mr=#1}{#2}
}
\providecommand{\href}[2]{#2}

\end{document}